\documentclass[10pt,a4paper]{amsart}
\usepackage{amsfonts}
\usepackage{amsthm}
\usepackage{amsmath}
\usepackage{mathtools}
\usepackage{amscd}
\usepackage[utf8]{inputenc}
\usepackage{t1enc}
\usepackage[mathscr]{eucal}
\usepackage{indentfirst}
\usepackage{graphicx}
\usepackage{graphics}
\usepackage{pict2e}
\usepackage{epic}
\usepackage{float}
\usepackage{MnSymbol}
\usepackage{multirow}
\usepackage{shuffle} 
\usepackage[table,xcdraw]{xcolor} 
\usepackage{hyperref}
\numberwithin{equation}{section}
\usepackage[margin=2.9cm]{geometry}
\usepackage{epstopdf} 
 
\allowdisplaybreaks

\newmuskip\pFqmuskip
\newcommand*\pFq[6][8]{%
	\begingroup 
	\pFqmuskip=#1mu\relax
	\mathchardef\normalcomma=\mathcode`,
	\mathcode`\,=\string"8000
	\begingroup\lccode`\~=`\,
	\lowercase{\endgroup\let~}\pFqcomma
	{}_{#2}F_{#3}{\left[\genfrac..{0pt}{}{#4}{#5};#6\right]}%
	\endgroup
}
\newcommand{\pFqcomma}{{\normalcomma}\mskip\pFqmuskip}

\theoremstyle{plain}
\newtheorem{theorem}{Theorem}[section]
\newtheorem{lemma}[theorem]{Lemma}
\newtheorem{corollary}[theorem]{Corollary}
\newtheorem{proposition}[theorem]{Proposition}

\theoremstyle{definition}

\newtheorem{problem}[theorem]{Problem}

\newcommand{\CMZV}[2]{\textsf{CMZV}^{#1}_{#2}}

\DeclareMathOperator{\Li}{Li}

\newtheorem{romexample}{Example}

\begin{document}
	
	\title{Multiple zeta values, WZ-pairs and infinite sums computations}
	
	\author[Kam Cheong Au]{Kam Cheong Au}
	
	\address{Rheinische Friedrich-Wilhelms-Universität Bonn \\ Mathematical Institute \\ 53115 Bonn, Germany} 
	
	\email{s6kmauuu@uni-bonn.de}
	

	\begin{abstract} We combine the powerful method of Wilf-Zeilberger pairs with systematic theory of multiple zeta values to prove a large number of series identities due to Z.W. Sun, many of them have been long standing conjectures.
	\end{abstract}
	
	\maketitle

	The combinatorialist Z.W. Sun conjectured (\cite{sun2013products}, \cite{sun2021book}, \cite{sun2010conjectures}, \cite{sun2022conjectures}) many infinite series from their $p$-adic analogous, such as
	$$\sum _{n=1}^{\infty } \frac{1}{n^2 \binom{2 n}{n}} (8 \sum _{j=1}^n \frac{(-1)^j}{j^4}+\frac{(-1)^n}{n^4})=-\frac{22 \zeta (3)^2}{15}-\frac{97 \pi ^6}{34020}$$
	$$\sum _{n=0}^{\infty } \frac{\binom{2 n}{n} }{(2 n+1)^2 (-16)^n} (40 \sum _{j=0}^{n} \frac{(-1)^j}{(2 j+1)^3}-\frac{47 (-1)^n}{(2 n+1)^3})=-\frac{85 \pi ^5}{3456}$$
	$$\sum_{n\geq 1} \frac{(28n^2-18n+3)(-64)^n}{n^5\binom{2n}{n}^4 \binom{3n}{n}} = -14\zeta(3)$$
	some of them are notoriously difficult to prove. In this article, we put forward a framework that can establish a large number of such conjectures: the main tools used are colored multiple zeta values (CMZVs) and Wilf-Zeilberger (WZ) pairs. 
	
	The technique of WZ-pair, developed in 1970s, is known for its ability to prove many identities. At first it was used to prove combinatorial identities (\cite{AequalsB}), later was generalized to infinite series (\cite{mohammed2004markov}, \cite{mohammed2005infinite}), able to prove series identities that otherwise quite hard to do (e.g. \cite{pilehrood2011bivariate}, \cite{pilehrood2008generating}, \cite{pilehrood2010series}, \cite{pilehrood2008simultaneous}, \cite{guillera2008hypergeometric}, \cite{guillera2010new}), with advantage that such proofs are all verifiable. \par
	
	The theory of MZVs, starting from 1980s, is now mature enough (e.g. \cite{gil2017multiple}, \cite{hoffman1997algebra}, \cite{deligne2010groupe}, \cite{zhao2016multiple}) to be applied back to classical problems like infinite sum, such approaches are found in, for instance, \cite{ablinger2014iterated}, \cite{ablinger2017discovering}, \cite{ablinger2019proving}, \cite{au2020evaluation}, \cite{au2022iterated}, \cite{xu2022sun}. Its techniques and applicable range are both very different from WZ method. \\[0.04in]
	
	Many of the "closed-form" conjectures of hypergeometric summation conform to the following shape (see Sun's articles \cite{sun2022conjectures}, \cite{sun2023new} for an avalanche of examples):
	\begin{equation}\label{overview}\sum_{n\geq 0} r^n\frac{(a_1)_n\cdots (a_p)_n}{(b_1)_n\cdots (b_p)_n} \times \left( \text{polynomials of harmonic numbers}\right) = \frac{\text{CMZVs of certain level}}{\pi^c}\end{equation}
	here "harmonic numbers" include not only the classical example $H_n^{(r)} = \sum_{i\leq n} i^{-r}$, but also others like $\sum_{i=1}^n \frac{(-1)^i}{i^4}$ and $\sum _{i=0}^{n} \frac{(-1)^i}{(2 i+1)^3}$, they can be viewed as partial sums of various versions of multiple polylogarithm (e.g. Hoffman's multiple $t$-values); $a_i,b_i$ are rational numbers in $(0,1]$, $r$ is an algebraic number, $c \in \mathbb{Z}^{\geq 0}$. 
	
	We briefly summarize currently known methods to attack these conjectures, they fall into four large, not necessarily mutually-exclusive, classes:
	\begin{itemize} \item Reduce them into CMZVs of some level, then invoke some datamine. (e,g. Ablinger \cite{ablinger2019proving} for Multiple Delinge Values, \cite{ablinger2017discovering}; later generalized by Au \cite{au2022iterated} to other levels)
		
		\item Functional equations of polylogarithm. (Charlton, Gangl, Xu, Zhao \cite{charlton2023two} is a typical manifestation of polylogarithm in such sums; Xu, Zhao \cite{xu2022sun}, \cite{xu2022note})
		
		\item WZ-method or its close relatives. They are used by Mohammed \& Pilehrood (e.g. \cite{mohammed2004markov}, \cite{mohammed2005infinite}, \cite{pilehrood2008generating}, \cite{pilehrood2008simultaneous}) and Guillera (e.g. \cite{guillera2008hypergeometric}). 
		
		\item Extracting coefficients from hypergeometric series with parameters. (Chu \cite{chu2020alternating}, \cite{chuAperyseries}, \cite{chu2014accelerating}; Wei \cite{wei2022conjectural}, \cite{wei2023some}, \cite{wei2023some2})
	\end{itemize}
	
	All four approaches have radically different natures. Up till recently, one did not know the connections between these methods. The idea will be presented in this article is a synthesis of 1st, 3rd and 4th method, it allows us to prove a plethora of Sun's conjectures. \par
	
	Nonetheless, there are still many unsolved problems. For example, all examples in this manuscript will have $c=0$ in \ref{overview}, but there are known examples with $c=1,2,3,4$. These so-called "Ramanujan $1/\pi$-formulas" are surveyed (without the harmonic number part) in \cite{cohen2021rational} by Cohen and Guillera, see also \cite{zudilin2007ramanujan} and \cite{zudilin2005quadratic} of Zudilin, from which Sun draws inspiration to produce more conjectures. While our method can adapt to some extend when $c\geq 1$, it seems not yet mature enough. Even in $c=0$ case, things are not settled, for example, 
	$$\sum_{n=1}^\infty 4^{-3n} \frac{(1)_n^3}{n^3(1/2)_n^3} (21 n-8) \left(\frac{43 H_{n-1}^{(3)}}{8}+H_{2 n-1}^{(3)}\right) \stackrel{?}{=} \frac{711 \zeta (5)}{28}-\frac{29 \pi ^2 \zeta (3)}{14}$$ remains conjectural (\textit{cf.} Example \ref{ex5} below).

	~\\[0.05in]
	\textbf{Acknowledgements} The author would like to express his sincere gratitude to W. Zudilin on boarding my perspective on these topics; D. Zagier, H. Gangl, S. Charlton and K. Bönisch for valuable discussions; C. Xu and J.Q. Zhao for introducing these problems.

	\section{Preliminaries}
	\subsection{Colored multiple zeta values}
	Consider the series
	$$\Li_{s_1,\cdots,s_k}(a_1,\cdots,a_k) = \sum_{n_1>\cdots>n_k\geq 1}\frac{a_1^{n_1}\cdots a_k^{n_k}}{n_1^{s_1} \cdots n_k^{s_k}}$$
	known as \textit{multiple polylogarithm}. Here $k$ is called the \textit{depth} and $s_1+\cdots+s_k$ is called the \textit{weight}. 
	
	When $a_i$ are $N$-th roots of unity, $s_i$ are positive integers and $(a_i, s_i) \neq (1,1)$, $\Li_{s_1,\cdots,s_k}(a_1,\cdots,a_k)$ is called a colored multiple zeta values (CMZV) of weight $s_1+\cdots+s_k$ and \textit{level} $N$. We say a complex number is CMZV of weight $w$ and level $N$ if it's in $\mathbb{Q}(e^{2\pi i /N})$-span of such numbers\footnote{in other literature, sometimes one takes $\mathbb{Q}$-span instead of $\mathbb{Q}(e^{2\pi i /N})$-span}, which we denote by $\CMZV{N}{w}$. The special case when all $a_i = 1$ is the well-known \textit{multiple zeta function}.
	~\\[0.01in]
	
	The motivic dimension of $\CMZV{N}{w}$ for small $N$ is known \cite{delignegroupes}. For small $N$ and $w$, we have an explicit database on expressing each such CMZVs into a linear combination of elements, whose number equals the motivic dimension. For $N=1, 2$, this is classical work of \cite{ihara2006derivation}, \cite{hoffman1997algebra}, Zhao (\cite{zhao2016multiple}, \cite{ZhaoStandard},\cite{zhao2008multiple}) made important progress on general $N$. Reaching the motivic dimension for general $N$ is still open, currently it's only known for $N\leq 8$ (\cite{au2022iterated}). Due to limited space, this article has the occasion to use $N=1,2,4,6$ only.
	
	For level $N=1,2$, an extensive database is the MZV datamine \cite{MZVdatamine}. For higher level, the only known (provable) reduction database\footnote{there are empirical databases which are based on very high precision numerical relations, for example, \cite{henn2017evaluating} for level 6 weight $\leq 6$ CMZVs} seems to be author's Mathematica package\footnote{which can be downloaded at https://www.researchgate.net/publication/357601353}. We will use such an explicit reduction into (a set of arbitrarily chosen) basis element\footnote{we remark that the word "basis" is misnomer, it's not known whether they're truly linear independent, but for explicit calculation, this does not matter.} for explicit calculation. 
	~\\[0.01in]
	
	We denote $G$ to be the Catalan constant. When $D$ is a fundamental discriminant, let $L_D(s) = \sum_{n\geq 1} \left(\frac{D}{n}\right) n^{-s}$ here $\left(\frac{D}{n}\right)$ is Jacobi symbol. When $D = -4$, we let $L_{-4}(s) = \beta(s)$, the Dirichlet beta function. It can be shown that $L_D(s)$ is a CMZV of level $|D|$ and weight $s$. 
	
	\subsection{Harmonic numbers and Euler sums}
	Fix a positive integer $N$, $0<\gamma_i \leq 1$ , we denote
	$$H(\vec{s},\vec{\gamma},N) = \sum_{\substack{n_1+\gamma_1 > n_2+\gamma_2 > \cdots > n_k+\gamma_k > 0 \\ n_1<N}} \frac{1}{(n_ 1+\gamma_1)^{s_ 1}\cdots (n_k+\gamma_k)^{s_k}}$$
	here $n_i$ in sum need to be integers; with $\vec{s} = (s_1,\cdots,s_n)\in 
	\mathbb{N}^n,\vec{\gamma} = (\gamma_1,\cdots,\gamma_k)$. 
	
	We write $\{a\}_n$ as $n$-component constant vector with entry $a$: $(a,a,\cdots,a)$. The special cases when $\vec{s}=\{s\}_n$ and $\vec{\gamma} = \{\gamma\}_n$ is denoted by $\hbar_N^{(s)}(\gamma)$, i.e. $$\hbar_N^{(s)}(\gamma) = \sum_{0\leq n\leq N-1} \frac{1}{(\gamma+n)^s}$$ This notation, when $\gamma=1$, is the familiar harmonic number $\hbar_N^{(s)}(1) = H_N^{(s)} = 1+1/2^s+\cdots+1/N^s$. 
	
	$\hbar_N^{(s)}(\gamma)$ is the most important type of harmonic number that features in our example, but many CMZVs-theoretical statement in next subsection also holds for $H(\vec{s},\vec{\gamma},N)$. \par
	
	It follows from Newton's identity of symmetric function that $H(\{s\}_k, \{\gamma\}_k, N)$ can be written as a (weighted) polynomial in terms of $\hbar_N^{is}(\gamma)$ for various $i$. For example:
	$$\begin{aligned}\label{bellpoly}
		H(\{s\}_1, \{\gamma\}_1, N) &= \hbar_N^{(s)}(\gamma) \\
		H(\{s\}_2, \{\gamma\}_2, N) &= \frac{1}{2}(\hbar_N^{(s)}(\gamma)^2 - \hbar_N^{(2s)}(\gamma)) \\
		H(\{s\}_3, \{\gamma\}_3, N) &= \frac{1}{6}(\hbar_N^{(s)}(\gamma)^3 - 3\hbar_N^{(s)}(1)\hbar_N^{(2s)}(\gamma) + 2\hbar_N^{(3s)}(\gamma))
	\end{aligned}$$
	the general formula is given by complete Bell polynomials \cite[Chapter~3]{comtet2012advanced}. We can define a more general harmonic sum than $H(\vec{s},\vec{\gamma},N)$:
	$$H(\vec{s},\vec{\gamma},\vec{a},N) = \sum_{\substack{n_1+\gamma_1 > n_2+\gamma_2 > \cdots > n_k+\gamma_k > 0 \\ n_1<N}} \frac{e(a_1 n_1)\cdots e(a_k n_k)}{(n_ 1+\gamma_1)^{s_ 1}\cdots (n_k+\gamma_k)^{s_k}}$$
	here $\vec{a} = (a_1,\cdots,a_k)$ and $e(x) = e^{2\pi i x}$. 
	
	\begin{theorem}\label{CMZVsum1}
		Let $0<\gamma_0\leq 1$, $K_i$ positive integers, $\vec{\gamma_i}$ and $\vec{a_i}$ rational vectors such that each component is $0< \cdot \leq 1$. Let $N$ be an integer such that all $N\vec{\gamma_i}, N\vec{a_i}, N\gamma_0, Na_0, N/K_i$ are in $\mathbb{Z}$. The following series, regularized (in sense defined below), is a level $N$ (possibly inhomogeneous) CMZV
		$$\sum_{n_0>0} \frac{e^{2\pi i  a_0 n_0}}{(n_0+\gamma_0)^{s_0}} H(\vec{s_1},\vec{\gamma_1},\vec{a_1},K_1 n_0) \cdots H(\vec{s_m},\vec{\gamma_m},\vec{a_m},K_m n_0) \in \sum_{1\leq i\leq w} \CMZV{N}{i}$$
		where $w = |\vec{s_1}| + \cdots + |\vec{s_i}| + s$. If $\gamma_0 = 1/N$, then it is a homogeneous CMZV of weight $w$. 
	\end{theorem}
	\begin{proof}
		We put the technical proof in appendix. 
	\end{proof}
	
	The proof also constitutes an algorithm, which is implemented in author's Mathematica package\footnote{this can be downloaded at https://www.researchgate.net/publication/357601353} as \textsf{MZSum}. Our ability to calculate these sums explicitly is important for \textit{effective proofs} of series evaluation in the main sections. Above theorem encompasses nearly all variations of Euler sums in literature that are related to CMZV (e.g. \cite{flajolet}, \cite{zhou2023sun}, \cite{zhou2022hyper}, \cite{xu2023dirichlet}).  \par
	
	We explain in what sense we are regularizing above sum. Let $M$ be a large integer, we will see in the proof that there exists $c_i \in \mathbb{C}$, \begin{multline*}\sum_{M\geq n_0>0} \frac{e^{2\pi i  a_0 n_0}}{(n_0+\gamma_0)^{s_0}} H(\vec{s_1},\vec{\gamma_1},\vec{a_1},K_1 n_0) \cdots H(\vec{s_m},\vec{\gamma_m},\vec{a_m},K_m n_0) = \\ \sum_{i>0} c_i (\log M + \gamma)^i + c_0 + o(1) \qquad M\to \infty\end{multline*}
	here $\gamma \approx 0.577$ is Euler's constant and only finitely many $c_i$ are non-zero. We define the regularized value to be $c_0$. 
	
	The technique of regularization is still crucial if we are only interested in convergent sum. For example, the convergent $\sum_{n\geq 0} \frac{H_n^3}{(4n+1)(4n+3)}$ is a CMZV of level 4, since we can split it (by partial fraction) into two divergent sums, and their regularized values are such CMZV. For the theory behind this, see \cite{racinet2002doubles}, \cite[Chapter~13]{zhao2008multiple}. 
	~\\[0.07in]

	\subsection{Series that gives CMZVs}
	Now we give examples that are immediate to our application. 
	Let $(a)_n = \Gamma(a+n)/\Gamma(a)$ be pochhammer symbol, for $0<\gamma\leq 1$, we have the series expansion (around $a=0$):
	\begin{equation}\label{pochhammerexpansion}(\gamma+a)_n = (\gamma)_n \left(1+ \sum_{k\geq 1} H(\{1\}_k, \{\gamma\}_k, n) a^k \right)\end{equation}
	To see this, write $$\frac{(\gamma+a)_n}{(\gamma)_n} = \prod_{i=0}^{n-1} \frac{\gamma+a+i}{\gamma+i} = \prod_{i=0}^{n-1} (1+\frac{a}{\gamma+i})$$
	the equation \ref{pochhammerexpansion} follows directly from definition of $H(\{1\}_k, \{\gamma\}_k,n)$. 
	
	Consider the function of $(a,b,c,d)$ \begin{equation}\label{ex0}\sum_{k\geq 0}\frac{(a+1)_k (b+1)_k}{(c+k+1) (d+k+1) (c+1)_k (d+1)_k}\end{equation}
	it's analytic near $(a,b,c,d)=(0,0,0,0)$, we wish to find its power series expansion, there seems no simple general formula for them, but we can make a qualitative statement: coefficient of $a^ib^jc^kd^l$ is a CMZV of level 1 with weight $i+j+k+l+2$, this follows from equation (\ref{pochhammerexpansion}). Using \ref{bellpoly}, the summand of each series can be expressed in terms of $H_n^{(r)}$ only. For example, coefficient of $ac^2d$ above is
	$$\sum_{k\geq 0}(-\frac{\left(H_k\right){}^2 H_k^{(2)}}{2 (k+1)^2}-\frac{H_k H_k^{(2)}}{2 (k+1)^3}-\frac{\left(H_k\right){}^4}{2 (k+1)^2}-\frac{3 \left(H_k\right){}^3}{2 (k+1)^3}-\frac{2 \left(H_k\right){}^2}{(k+1)^4}-\frac{H_k}{(k+1)^5})$$
	by converting to our CMZV basis, it is $-\frac{5 \zeta (3)^2}{2}-\frac{3 \pi ^6}{140}$. More generally, expanding $\ref{ex0}$ at point $a=b=c=d = -1+r/s$ with $0<r<s$, the coefficients are level $s$ CMZV. 
	
	Actually, there are much more series which produce CMZV, for example
	\begin{proposition}\label{binomCMZVs}
		$$\sum_{k\geq 0} \frac{1}{(k+1/2)^s} \frac{(1)_n}{(1/2)_n}\hbar_n^{(r_1)}(1)\cdots \hbar_n^{(r_i)}(1) \hbar_n^{(s_1)}(1/2)\cdots \hbar_n^{(s_k)}(1/2)$$
		is CMZV of level 4, of weight $s+\sum r_i + \sum s_i$. 
		There exists similar assertion for $(1/2)_n/(1)_n$. 
	\end{proposition}
	\begin{proof}
		The proof is non-trivial and involves much about general CMZV theory. See \cite{au2020evaluation}, \cite{au2022iterated}, \cite{davydychev2004binomial}, \cite{kalmykovBinomial}, \cite{xu2022ap} to get an idea on how to handle such harmonic sum twisted with binomial coefficient. This particular assertion, as well as more general ones, will be proved in an upcoming article of author. 
	\end{proof}
	
	Since it uses not-yet published result, we refrain from using above proposition, it only appears in Example \ref{ex9} in this article.
	
	\section{Wilf-Zeilberger pairs}
	For introduction and terminology, see \cite{AequalsB}, \cite{mohammed2005infinite}. Let $F(n,k), G(n,k)$ be a WZ-pair, that is, they satisfy the relation
	$$F(n+1,k)-F(n,k) = G(n,k+1)-G(n,k)$$
	
	\begin{proposition}[\cite{mohammed2005infinite}]\label{WZsumformula} Let $F,G$ as above, suppose \begin{itemize}\item $\lim_{k\to\infty} G(n,k)=0$ for each $n\geq 0$ 
			\item $\sum_{k\geq 0} F(0,k)$ converges
			\item $\sum_{n\geq 0} G(n,0)$ converges
		\end{itemize}
		then $\lim_{n\to\infty} \sum_{k\geq 0} F(n,k)$ also exists and we have $$\sum_{k\geq 0} F(0,k) = \sum_{n\geq 0} G(n,0) + \lim_{n\to\infty} \sum_{k\geq 0} F(n,k)$$
	\end{proposition}
	\begin{proof}
		We apply $\sum_{n= 0}^{N-1} \sum_{k=0}^{K-1}$ on both sides of $F(n+1,k)-F(n,k) = G(n,k+1)-G(n,k)$, giving
		$$\sum_{k=0}^{K-1} (F(N,k)-F(0,k)) = \sum_{n= 0}^{N-1} (G(n,K) - G(n,0))$$
		by our assumption, we can let $K\to \infty$ while fixing $N$, to obtain
		$$\sum_{k\geq 0} F(N,k) - \sum_{k\geq 0} F(0,k) = -\sum_{n= 0}^{N-1} G(n,0)$$
		then letting $N\to\infty$ proves the claim.
	\end{proof}
	
	Given a proper hypergeometric $F$, Gosper's algorithm can be used to determine whether such proper hypergeometric $G$ exist. But there is an empirical method to generate many $F$'s for which $G$ does exist:
	
	\begin{problem}[Hypergeometric summation formula induces WZ-candidate]\label{summation_induces_WZ}
		Let $f(a,b,...;k)$ be a proper hypergeometric term in $a,b,...,k$. If $$\sum_{k\geq 0} f(a,b,...;k)$$ exists and is independent of $a,b,...$, then is it true that, for any integers $A,B,...,K$,
		$$F(n,k) = f(a+An,b+Bn,...;k+Kn)$$
		has a proper hypergeometric WZ-mate $G$?
	\end{problem}
	
	The answer seems to be affirmative for all cases we have computed. Gessel (\cite{gessel1995finding}) already used this empirical observation to generate a large amount of finite hypergeometric identities. We give an example of above observation, consider the Gauss $_2F_1$ formula:
	$$\pFq{2}{1}{a \quad b}{c}{1} = \sum_{k\geq 0} \frac{(a)_k(b)_k}{(1)_k(c)_k} = \frac{\Gamma(c-a-b)\Gamma(c)}{\Gamma(c-a)\Gamma(c-b)}$$
	divide the gamma factor of RHS into LHS summand, we arrive:
	$$\sum_{k\geq 0} \frac{\Gamma (c-a) \Gamma (a+k) \Gamma (c-b) \Gamma (b+k)}{\Gamma (a) \Gamma (b) \Gamma (k+1) \Gamma (c+k) \Gamma (-a-b+c)} = 1$$
	If we denote the summand by $\textsf{Gauss2F1}(a,b,c;k)$, then the above observation would predict $$F(n,k) = \textsf{Gauss2F1}(a - d - n, b - d - n, 1 + c - d; 1 + d + k + n)$$
	has a WZ-mate $G$ that is also proper hypergeometric, indeed, one verifies
	\begin{multline*}G(n,k) =F(n,k) \times (a b-a c-a d-a k-2 a n-2 a-b c-b d-b k-2 b n-2 b+c^2+c d+c k+3 c n+3 c \\ +d^2+d k+3 d n+3 d+2 k n+2 k+3 n^2+6 n+3) \div ((a+b-c-d-2 n-2) (a+b-c-d-2 n-1)) \end{multline*}
	satisfies $F(n+1,k)-F(n,k) = G(n,k+1)-G(n,k)$.
	
	Apart from $\textsf{Gauss2F1}(a,b,c;k)$ defined above, we shall also define
	\begin{align*}&\textsf{Dixon3F2}(a,b,c;k)  \\ &= \frac{2 \Gamma (a-b+1) \Gamma (a-c+1) \Gamma (2 a+k) \Gamma (b+k) \Gamma (c+k) \Gamma (2 a-b-c+1)}{\Gamma (a) \Gamma (b) \Gamma (c) \Gamma (k+1) \Gamma (a-b-c+1) \Gamma (2 a-b+k+1) \Gamma (2 a-c+k+1)}\end{align*}
	\begin{align*}&\textsf{Dougall5F4}(a,b,c,d;k) \\ &= \frac{(a+2 k) \Gamma (a+k) \Gamma (b+k) \Gamma (c+k) \Gamma (d+k) \Gamma (a-b-c+1) \Gamma (a-b-d+1) \Gamma (a-c-d+1)}{\Gamma (b) \Gamma (c) \Gamma (d) \Gamma (k+1) \Gamma (a-b+k+1) \Gamma (a-c+k+1) \Gamma (a-d+k+1) \Gamma (a-b-c-d+1)}\end{align*}
	as well as
	\begin{align*}&\textsf{Watson3F2}(a,b,c;k) \\ &= \frac{2^{-2 a-2 b+2 c+1} \Gamma \left(-a+c+\frac{1}{2}\right) \Gamma (2 a+k) \Gamma \left(-b+c+\frac{1}{2}\right) \Gamma (2 b+k) \Gamma (c+k)}{\Gamma (a) \Gamma (b) \Gamma (k+1) \Gamma (2 c+k) \Gamma \left(-a-b+c+\frac{1}{2}\right) \Gamma \left(a+b+k+\frac{1}{2}\right)}
	\end{align*}
	in all cases their sums over $k\geq 0$ are $1$:\footnote{subjected to restriction on $a,b,c,d$ such that the sum converges.}
	$$\begin{aligned}\sum_{k\geq 0} \textsf{Gauss2F1}(a,b,c;k) &= 1 \\
		\sum_{k\geq 0} \textsf{Dixon3F2}(a,b,c;k) &= 1 \\
		\sum_{k\geq 0} \textsf{Watson3F2}(a,b,c;k) &= 1 \\
		\sum_{k\geq 0} \textsf{Dougall5F4}(a,b,c,d;k) &= 1\end{aligned}$$
	these follow by using the corresponding $_pF_q$ summation formulas. 
	
	\par There are many computer algebra packages to find WZ-pairs, such as the Maple package in \cite{pilehrood2011bivariate}, \cite{pilehrood2008simultaneous}. Mathematica is used for all calculations, both CMZVs and WZ-pairs, in this manuscript. The WZ-pairs package we use comes from \href{https://www3.risc.jku.at/research/combinat/software/ergosum/RISC/fastZeil.html}{RISC} (i.e. Gosper's algorithm). 
	
	\subsection{The boundary term \texorpdfstring{$\lim_{n\to\infty} \sum_{k\geq 0} F(n,k)$}{TEXT}}
	
	In all examples of this article, $\lim_{k\to\infty} G(n,k)=0$ is always satisfied. For examples in first two sections, $\lim_{n\to\infty} \sum_{k\geq 0} F(n,k) = 0$ is also satisfied, so we have, according to \ref{WZsumformula}, an elegant equality $\sum_{k\geq 0} F(0,k) = \sum_{n\geq 0} G(n,0)$. We will also give some sporadic examples (Examples \ref{ex9}, \ref{ex10}) with non-zero $\lim_{n\to\infty} \sum_{k\geq 0} F(n,k)$, in these cases, computing the limit requires some tricks in asymptotic analysis. \par
	
	We state a sufficient condition for $\lim_{n\to\infty} \sum_{k\geq 0} F(n,k) = 0$ that will be convenient to us. We let $\tilde{\Gamma}(x) = \Gamma(a+x)$ for some unspecified complex constant $a$. different occurrences of $\tilde{\Gamma}$ need not have the same $a$. 
	\begin{proposition}\label{vanishingboundary}
		Let $A_i,B_i,C_i,D_i>0$. Consider the limit
		$$\lim_{n\to\infty} \tilde{\Gamma}\begin{pmatrix}A_1n & \cdots &A_k n \\ B_1n & \cdots & B_m n \end{pmatrix}\sum_{k\geq 0} \tilde{\Gamma}\begin{pmatrix}k+C_1n & \cdots & k+C_ln \\ k+D_1n & \cdots & k+ D_l n\end{pmatrix}$$ If $\sum_i A_i + \sum_i C_i = \sum_i B_i + \sum_i D_i$ and $$\frac{A_1^{A_1}\cdots A_k^{A_k}}{B_1^{B_1}\cdots B_m^{B_m}} \frac{(x+C_1)^{x+C_1} \cdots (x+C_l)^{x+C_l}}{(x+D_1)^{x+D_1} \cdots (x+D_l)^{x+D_l}} < 1 \qquad \forall x\in \mathbb{R}^{\geq 0}$$
		holds, then above limit is $0$. 
	\end{proposition}
	\begin{proof}
		Write $f(n)\ll_P g(n)$ if $f(n)/g(n)$ is at most polynomial growth in $n$ as $n\to \infty$. We note that for $E_i, F_i>0$ with $\sum E_i = \sum F_i$, $$\tilde{\Gamma}\begin{pmatrix}E_1n & \cdots &E_k n \\ F_1n & \cdots & F_m n \end{pmatrix} \ll_P \left(\frac{E_1^{E_1}\cdots E_k^{E_k}}{F_1^{F_1}\cdots F_m^{F_m}}\right)^n$$
		with the $\ll_P$ uniform in $E_i, F_i>0$. This assertion follows easily from Stirling's formula. For our original series, we split $\sum_{k\geq 0}$ into $\sum_{t\geq 0} \sum_{tn \leq k < (t+1)n}$, there exists $0\leq \theta_{t,n}\leq 1$ such that
		$$\sum_{tn\leq k < (t+1)n} \tilde{\Gamma}\begin{pmatrix}k+C_1n & \cdots & k+C_ln \\ k+D_1n & \cdots & k+ D_l n\end{pmatrix} = n \tilde{\Gamma}\begin{pmatrix}(t+\theta_{t,n}+C_1)n & \cdots & (t+\theta_{t,n}+C_l)n \\ (t+\theta_{t,n}+D_1)n & \cdots & (t+\theta_{t,n}+D_l) n\end{pmatrix}$$
		therefore the original limit is $$\ll_P \sum_{t\geq 0} \left(\frac{A_1^{A_1}\cdots A_k^{A_k}}{B_1^{B_1}\cdots B_m^{B_m}} \frac{(t+\theta_{t,n}+C_1)^{t+\theta_{t,n}+C_1} \cdots (t+\theta_{t,n}+C_l)^{t+\theta_{t,n}+C_l}}{(t+\theta_{t,n}+D_1)^{t+\theta_{t,n}+D_1} \cdots (t+\theta_{t,n}+D_l)^{t+\theta_{t,n}+D_l}} \right)^n$$
		Each number in parenthesis is $<1$, monotone convergence theorem implies it is $0$ as $n\to \infty$. 
	\end{proof}
	
	As an immediate application, above proposition implies $$\lim_{n\to \infty} \sum_{k\geq 0}\textsf{Gauss2F1}(a,b,c+n;d+n+k) = 0$$
	for any $a,b,c,d\in \mathbb{C}$. Indeed, we only need to check
	whether $$\frac{1^1 1^1}{1^1} \frac{(1+x)^{1+x} (1+x)^{1+x}}{(1+x)^{1+x} (2+x)^{2+x}} < 1$$
	for all non-negative $x$, this is evident.
	
	\section{Introductory examples}
	
	Now we combine techniques from previous sections to derive many series summation conjectures.

	The author believes understanding the conceptual idea of this technique is more important than the actual computational details, which are tedious and not quite do-able by hand even in simplest cases, therefore, most such computations are suppressed and only final result is presented. \par
	Readers interested in these computation details can consult accompanying Mathematica notebooks. 
	\begin{romexample}\label{ex1}For $\Re(a+b-c-d)<1$, we have
		$$\sum_{k\geq 0}\frac{(a+1)_k (b+1)_k}{(c+1)_{k+1} (d+1)_{k+1}} = \sum_{n\geq 1} \frac{(a+1)_{n-1} (b+1)_{n-1} (-a+c+1)_{n-1} (-b+c+1)_{n-1} \times P}{(c+1)_{2 n} (d+1)_n (-a-b+c+d+1)_n} $$
		here $P =  a b-a c-a n-b c-b n+c^2+c d+3 c n+2 d n+3 n^2$. 
	\end{romexample}
	\begin{proof}
		Let $F(n,k) = \textsf{Gauss2F1}(a - d, b - d, 1 + c - d + n; 1 + n + k + d)$, one verifies it has a proper hypergeometric WZ-mate $G(n,k)$ (as predicted by Problem \ref{summation_induces_WZ}). The sum $\sum_{k\geq 0} F(0,k)$ converges if $\Re(a+b-c-d)<1$; $\sum_{n\geq 0} G(n,0)$ always converges regardless of $a,b,c,d$; the boundary term $\lim_{n\to\infty} \sum_{k\geq 0} F(n,k)$ vanishes by Proposition \ref{vanishingboundary}. So Proposition \ref{WZsumformula} implies $\sum_{k\geq 0} F(0,k) = \sum_{n\geq 1} G(n-1,0)$, which is exactly the above formula after dividing by gamma factors that are independent of $n$ and $k$.
	\end{proof}
	
	We will use the template of the proof so often that it is convenient to call above displayed formula \textit{the summation formula induced by} $F(n,k) = \textsf{Gauss2F1}(a - d, b - d, 1 + c - d + n; 1 + n + k + d)$.
	
	Expanding Example \ref{ex1} near $(a,b,c,d)=(0,0,0,0)$, coefficients of $a^i b^j c^k d^l$ of LHS are CMZVs of level 1 and weight $i+j+k+l+2$. For RHS, we extract the coefficient using  Formula \ref{pochhammerexpansion}, we obtain equalities of the form 
	$$\text{level 1 homogeneous MZV} = \sum_{n\geq 1} \frac{(1)_n^2}{(1)_{2n}} \mathbb{Q}[n,\frac{1}{n},H^{(r)}_n,H^{(r)}_{2n}]\qquad r = 1,2,\cdots$$
	we call $\frac{(1)_n^2}{(1)_{2n}} = 4^{-n} \frac{(1)_n}{(\frac{1}{2})_n}$ the hypergeometric part of the formula, this is $\binom{2n}{n}^{-1}$ in terms of binomial coefficient. 
	We obtain, for instances, the following identities by comparing coefficient of $1,a,c,a^2$ respectively:
	$$\begin{aligned}\frac{\pi ^2}{6} &= \sum_{n\geq 1} 4^{-n} \frac{(1)_n}{(\frac{1}{2})_n} \frac{3}{n^2} \\
		\zeta(3) &= \sum_{n\geq 1} 4^{-n} \frac{(1)_n}{(\frac{1}{2})_n} (\frac{3 H_n}{n^2}-\frac{1}{n^3}) \\
		-2\zeta(3) &= \sum_{n\geq 1} 4^{-n} \frac{(1)_n}{(\frac{1}{2})_n} (\frac{3 H_n}{n^2}-\frac{3 H_{2 n}}{n^2}-\frac{3}{n^3}) \\
		\frac{\pi^4}{90} &= \sum_{n\geq 1} 4^{-n} \frac{(1)_n}{(\frac{1}{2})_n} (-\frac{H_n}{n^3}+\frac{3 \left(H_n\right){}^2}{2 n^2}-\frac{3 H_n^{(2)}}{2 n^2}+\frac{3}{n^4})\end{aligned} $$
	
	See \cite{chuAperyseries} for more similar examples involving hypergeometric part $4^{-n} \frac{(1)_n}{(\frac{1}{2})_n}$. For a method with CMZV nature applied to above example, see \cite{ablinger2017discovering} and \cite{ablinger2019proving}, which uses level $6$ CMZVs. \cite{bailey2006experimental} is an empirical approach to these identities.

	\begin{romexample}\label{ex2}
		For $a,b,c,d,e$ near $0$, we have
		\begin{multline*}\sum _{k\geq 0} \frac{(2+e+2 k) (1+a)_k (1+b)_k (1+c)_k (1+d)_k}{(1-a+e)_{k+1} (1-b+e)_{k+1} (1-c+e)_{k+1} (1-d+e)_{k+1}}= \\ \sum_{n\geq 1}\frac{\splitfrac{(-1)^n P (1-a-b+e)_{n-1} (1-a-c+e)_{n-1} (1-b-c+e)_{n-1}}{ (1-a-d+e)_{n-1} (1-b-d+e)_{n-1} (1-c-d+e)_{n-1} }}{(1-a+e)_n (1-b+e)_n (1-c+e)_n (1-d+e)_n (1-a-b-c-d+2 e)_{2 n}} \end{multline*}
		here \small{$P = a^2 b+a b^2+a^2 c+2 a b c+b^2 c+a c^2+b c^2+a^2 d+2 a b d+b^2 d+2 a c d+2 b c d+c^2 d+a d^2+b d^2+c d^2-2 a^2 e-5 a b e-2 b^2 e-5 a c e-5 b c e-2 c^2 e-5 a d e-5 b d e-5 c d e-2 d^2 e+6 a e^2+6 b e^2+6 c e^2+6 d e^2-5 e^3-2 a^2 n-5 a b n-2 b^2 n-5 a c n-5 b c n-2 c^2 n-5 a d n-5 b d n-5 c d n-2 d^2 n+12 a e n+12 b e n+12 c e n+12 d e n-15 e^2 n+6 a n^2+6 b n^2+6 c n^2+6 d n^2-15 e n^2-5 n^3$}.
	\end{romexample}
	\begin{proof}
		This is the formula induced from $F(n,k) = \textsf{Dougall5F4}(2 a - e - n, a + b - e - n, a + c - e - n, a + d - e - n, 1 - a + e + k + n)$, the boundary term $\lim_{n\to\infty} \sum_{k\geq 0}F(n,k) = 0$ by Proposition \ref{vanishingboundary}. The formula is $\sum_{k\geq 0}F(0,k) = \sum_{n\geq 1}G(n-1,0)$.
	\end{proof}
	
	The above is equivalent to Theorem 10 in \cite{chu2014accelerating}. When expanding at $(a,b,c,d,e)=(0,0,0,0,0)$, LHS is homogeneous CMZV of level $1$, so we obtain equalities of the form
	$$\text{level 1 homogeneous MZV} = \sum_{n\geq 1} (-4)^{-n} \frac{(1)_n}{(1/2)_{n}} \mathbb{Q}[n,\frac{1}{n},H^{(r)}_n,H^{(r)}_{2n}]\qquad r = 1,2,\cdots$$
	
	In \cite{chu2020alternating}, this method of comparing coefficients has already been used on Example \ref{ex2} to generate many identities, such as $$\begin{aligned}\sum_{n\geq 1} \frac{(-1)^{n-1}}{\binom{2n}{n}} (\frac{10 H_n}{n^3} - \frac{3}{n^4}) &= \frac{\pi^4}{30} \\ \sum_{n\geq 1} \frac{(-1)^{n-1}}{\binom{2n}{n}} (\frac{4 H_n}{n^3} + \frac{H_{2n}}{n^3}) &= \frac{2\pi^4}{75} \\
		\sum_{n\geq 1} \frac{(-1)^{n-1}}{\binom{2n}{n}} (\frac{1}{n^6} + \frac{5 H_n^{(3)}}{n^3}) &= 2\zeta(3)^2 \end{aligned}$$
	they are obtained as certain linear combinations when comparing coefficients of above formula. For example, the last one is $3[a^3]-3[a^2b]-[abc]$, where $[..]$ means the equality obtained by comparing the corresponding coefficient. 
	
	Sun \cite{sun2022conjectures} recently conjectured $$\sum _{n=1}^{\infty } \frac{(-1)^{n-1}}{n^3 \binom{2 n}{n}}\left(H_{2 n-1}^{(2)}-\frac{123 H_{n-1}^{(2)}}{16}\right)=\frac{451 \zeta (5)}{40}-\frac{14 \pi ^2 \zeta (3)}{15}$$
	despite matching the form, this one \textit{does not follow} from above formula. To prove this, we need another formula.
	
	\begin{romexample}\label{ex3}
		For $a,b,c,d$ near $0$, we have
		\begin{multline*}\sum _{k\geq 0} \frac{2(1+2 b+2 d+2 k) (1+c)_k (1+2 a+d)_k (\frac{1}{2}+b+d)_k}{(1+d)_{k+1} (\frac{1}{2}+2 a-b+d)_{k+1} (1+2 a-c+2 d)_{k+1}} = \\ -\sum _{n\geq 0} \frac{(-1)^n P (1+2 a+d)_{n-1} (\frac{1}{2}+b+d)_n (1-c+d)_{n-1} (1+a-c+d)_{n-1} (\frac{1}{2}+2 a-b-c+d)_{n-1}}{ (1+d)_n (\frac{1}{2}+2 a-b+d)_n (\frac{1}{2}+a-b-c+d)_n (1+2 a-c+2 d)_{2 n}} \end{multline*}
		here $P = -2 a+8 a^2-4 a b+2 c-12 a c+4 b c+4 c^2-3 d+20 a d-6 b d-12 c d+10 d^2-3 n+20 a n-6 b n-12 c n+20 d n+10 n^2$.
	\end{romexample}
	\begin{proof}
		This is the formula induced from $F(n,k) = \textsf{Dixon3F2}(a, b + 1/2, c - d - n, 1 + k + n + d)$, the boundary term $\lim_{n\to\infty} \sum_{k\geq 0} F(n,k) = 0$ by Proposition \ref{vanishingboundary}. The formula is $\sum_{k\geq 0}F(0,k) = \sum_{n\geq 1}G(n-1,0)$.
	\end{proof}
	
	When expanding at $(a,b,c,d)=(0,0,0,0)$, LHS is inhomogeneous CMZV of level $2$,\footnote{actually they are \textit{homogeneous of level 1}, but showing this requires a deeper development of above method, which we postpone to a later article.} so we obtain equalities of the form
	$$\text{level 2 inhomogeneous MZV} = \sum_{n\geq 1} (-4)^{-n} \frac{(1)_n}{(1/2)_{n}} \mathbb{Q}[n,\frac{1}{n},H^{(r)}_n,H^{(r)}_{2n}]\qquad r = 1,2,\cdots$$
	
	\begin{corollary}
		$$\sum _{n=1}^{\infty } \frac{(-4)^{-n} (1)_n}{(1/2)_n n^3}\left(H_{2 n-1}^{(2)}-\frac{123 H_{n-1}^{(2)}}{16}\right)=\frac{451 \zeta (5)}{40}-\frac{14 \pi ^2 \zeta (3)}{15}$$
	\end{corollary}
	\begin{proof}
		If we write $a^ib^jc^kd^l$ as the equality obtained from above formula by comparing corresponding coefficient, the above is the following linear combination: $-\frac{18501 a^3}{1520}-\frac{28343 a^2 b}{3040}+\frac{7567 a^2 c}{1520}+\frac{1505 a^2 d}{608}+\frac{23173 a b^2}{3040}+\frac{2609 a b c}{760}+\frac{3943 a b d}{608}+\frac{2653 a c^2}{760}-\frac{161 a c d}{304}+\frac{2233 a d^2}{608}+\frac{6589 b^3}{380}+\frac{24497 b^2 c}{1520}+\frac{3977 b c^2}{760}+\frac{81 b c d}{32}-\frac{2212 c^3}{95}$.
	\end{proof}
	
	Due to the complicated coefficients involved in above example, we say a few words on how they are found. For each $a^ib^jc^kd^l$ with $i+j+k+l = 3$, this coefficient of RHS is of form $\sum_{n\geq 0} f(n)$, with $f: \mathbb{N} \to \mathbb{Q}$. After we computed all such $f$ corresponding to $N$ monomials of degree $3$, then we have $N$ such equations
	$$\sum_{n\geq 0} f_1(n) \qquad \sum_{n\geq 0} f_2(n)\qquad \cdots \qquad \sum_{n\geq 0}f_N(n)$$
	our concern is whether there exists $c_i\in \mathbb{Q}$ such that
	\begin{equation}\label{cisolve}\sum_{i=1}^N c_i f_i(n) = \frac{(-4)^{-n} (1)_n}{(1/2)_n n^3}\left(H_{2 n-1}^{(2)}-\frac{123 H_{n-1}^{(2)}}{16}\right)\end{equation}
	treating $c_i$ as unknowns, plugging in $n=1,2,\cdots$ to obtain as many as linear equations in $c_i$ as we desire, Gaussian elimination is then used to found out effectively whether this over-determined inhomogenous linear system over $\mathbb{Q}$ has solution. When such $c_i$'s indeed exist by solving the system, validity of equation \ref{cisolve} for general $n$ is straightforward to verify symbolically.
	
	We will subsequently employ the above strategy without further comment.\\[0.07in]
	
	We can also compare coefficient on Example \ref{ex1} at the point $(a,b,c,d)=(-1/2,-1/2,-1/2,-1/2)$. LHS will be homogeneous CMZVs of level 2, RHS will be infinite series with hypergeometric part $4^{-n} \frac{(1)_n^3}{(1/2)_n^3}$ and kind of harmonic numbers involved are $H_n^{(r)}$ and $H_{2n}^{(r)}$. In particular, we have
	
	\begin{corollary}
		$$\sum_{n\geq 1} \frac{4^{-n}(1)_n^3}{(1/2)_n^3} (\frac{3 n-1}{n^3}) = \frac{\pi^2}{2} $$
		$$\sum_{n=1}^\infty  \frac{4^{-n}(1)_n^3}{n^3(1/2)_n^3} (H_{2 n-1}^{(2)}-\frac{5 H_{n-1}^{(2)}}{4}) = \frac{\pi^4}{24}$$
	\end{corollary}
	
	First is originally proved in \cite{guillera2008hypergeometric}; second is also proved in \cite{wei2022conjectural}, using essentially same method and formula \ref{ex1}.

	\section{More examples}
	\begin{romexample}\label{ex4}For $a,b,c,d$ near $0$, we have
		\begin{multline}16\sum_{k\geq 0} \frac{(c+1)_k (2 a+d+\frac{1}{2})_k (b+d+\frac{1}{2})_k}{(d+\frac{1}{2})_{k+1} (2 a-b+d+\frac{1}{2})_{k+1} (2 a-c+2 d+1)_k} \\ = \sum_{n\geq 0} \frac{(-1)^n (2 a+d+\frac{1}{2})_n (b+d+\frac{1}{2})_n (-c+d+\frac{1}{2})_n (a-c+d+\frac{1}{2})_n (2 a-b-c+d+\frac{1}{2})_n \times P}{ (d+\frac{1}{2})_{n+1} (2 a-b+d+\frac{1}{2})_{n+1} (2 a-c+2 d+1)_{2 n+1} (a-b-c+d+\frac{1}{2})_{n+1}}\end{multline}
		where $P = 16 a^2-8 a b-24 a c+40 a d+40 a n+20 a+8 b c-12 b d-12 b n-6 b+8 c^2-24 c d-24 c n-12 c+20 d^2+40 d n+20 d+20 n^2+20 n+5$.
	\end{romexample}
	\begin{proof}
		This is the formula induced from $F(n,k) = \textsf{Dixon3F2}(a,b,c-d-n+\frac{1}{2},d+k+n+\frac{1}{2})$, the boundary term $\lim_{n\to\infty}\sum_{k\geq 0}F(n,k) = 0$ by Proposition \ref{vanishingboundary}. The formula is $\sum_{k\geq 0}F(0,k) = \sum_{n\geq 0}G(n,0)$.
	\end{proof}
	
	Expanding LHS at $(a,b,c,d) = (0,0,0,0)$, coefficients are homogeneous CMZVs of level 2. RHS's coefficients are series form
	$$\sum_{n\geq 0} (-\frac{1}{4})^n \frac{(1/2)_n}{(1)_n} \mathbb{Q}[\frac{1}{2n+1}, \hbar_n^{(r)}(1), \hbar_n^{(r)}(1/2)] \qquad r=1,2,\cdots$$ in terms of binomial coefficient, the hypergeometric part is $(-16)^{-n}\binom{2n}{n}$. The constant term gives $$\sum_{n\geq 0} (-\frac{1}{4})^n \frac{(1/2)_n}{(1)_n}\frac{5}{(2 n+1)^2} = \frac{\pi^2}{2}$$
	
	\begin{corollary}[Conjectures in \cite{sun2021book}, \cite{sun2010conjectures}]
		$$\begin{aligned}\sum_{n\geq 0} (-\frac{1}{4})^n \frac{(1/2)_n}{(1)_n} \frac{1}{(2 n+1)^2}  \left(5H_{2n+1} + \frac{12}{2n+1}\right) &= 14\zeta(3)\\
			\sum_{n\geq 0} (-\frac{1}{4})^n \frac{(1/2)_n}{(1)_n} \frac{1}{(2 n+1)^2}  \left(5 \sum _{j=0}^n \frac{1}{(2 j+1)^3}+\frac{1}{(2 n+1)^3}\right) &= \frac{\pi^2}{2}\zeta(3)\\
			\sum_{n\geq 0} (-\frac{1}{4})^n \frac{(1/2)_n}{(1)_n} \frac{1}{(2 n+1)^2}  \left(5 \sum _{j=0}^{n-1} \frac{1}{(2 j+1)^4}+\frac{1}{(2 n+1)^4}\right) &= \frac{7\pi^6}{7200}\end{aligned}$$
	\end{corollary}
	\begin{proof}
		They all follow from comparing coefficients of above formula.
	\end{proof}
	The first equality above was also proved in \cite{charlton2023two} via functional equations of polylogarithm. 
	
	\begin{romexample}\label{ex5}For $a,b,c,d$ near $0$, we have
		\begin{multline*}\sum _{k\geq 0} \frac{(1+a)_k (1+b)_k}{(1+c)_{k+1} (1+d)_{k+1}} \\  = \sum_{n\geq 1} \frac{P (1+a)_{n-1} (1+b)_{n-1} (1-a+c)_{n-1} (1-b+c)_{n-1} (1-a+d)_{n-1} (1-b+d)_{n-1}}{(1+c)_{2 n} (1+d)_{2 n} (1-a-b+c+d)_{2 n}}\end{multline*}
		here \small{$P = a^2 b^2-a^2 b c-a b^2 c+a b c^2-a^2 b d-a b^2 d+a c d+a^2 c d+b c d+3 a b c d+b^2 c d-c^2 d-2 a c^2 d-2 b c^2 d+c^3 d+a b d^2-c d^2-2 a c d^2-2 b c d^2+2 c^2 d^2+c d^3-2 a^2 b n-2 a b^2 n+2 a c n+a^2 c n+2 b c n+6 a b c n+b^2 c n-2 c^2 n-3 a c^2 n-3 b c^2 n+2 c^3 n+2 a d n+a^2 d n+2 b d n+6 a b d n+b^2 d n-6 c d n-11 a c d n-11 b c d n+10 c^2 d n-2 d^2 n-3 a d^2 n-3 b d^2 n+10 c d^2 n+2 d^3 n+4 a n^2+a^2 n^2+4 b n^2+8 a b n^2+b^2 n^2-8 c n^2-13 a c n^2-13 b c n^2+13 c^2 n^2-8 d n^2-13 a d n^2-13 b d n^2+29 c d n^2+13 d^2 n^2-8 n^3-14 a n^3-14 b n^3+28 c n^3+28 d n^3+21 n^4$}.
	\end{romexample}
	\begin{proof}
		This is the formula induced by $\textsf{Gauss2F1}(a-d-n,b-d-n,c-d+1,d+k+2 n+1)$, boundary term $\lim_{n\to\infty}\sum_{k\geq 0}F(n,k) = 0$ by Proposition \ref{vanishingboundary}. The formula is $\sum_{k\geq 0}F(0,k) = \sum_{n\geq 1}G(n-1,0)$.
	\end{proof}

	Compare coefficient on Example \ref{ex5} at the point $(a,b,c,d)=(0,0,0,0)$. LHS will be homogeneous CMZVs of level 1, RHS will be infinite series with hypergeometric part $4^{-3n} \frac{(1)_n^3}{(1/2)_n^3}$ and kind of harmonic numbers involved are $H_n^{(r)}$ and $H_{2n}^{(r)}$. Taking certain linear combinations of above formula, we obtain results such as
	
	\begin{corollary}
		$$\begin{aligned}\sum_{n\geq 1} 4^{-3n} \frac{(1)_n^3}{(1/2)_n^3} (\frac{21}{n^2}-\frac{8}{n^3}) &= \zeta(2)\\
			\sum_{n=1}^\infty 4^{-3n} \frac{(1)_n^3}{n^3(1/2)_n^3} (H_{2 n-1}^{(2)}-\frac{25 H_{n-1}^{(2)}}{8}) &= \frac{47\pi^4}{2880}\end{aligned}$$
	\end{corollary}
	
	The first one of above corollary are a famous WZ-type result due to Zeilberger \cite{zeilberger1993closed}, from which Sun \cite{sun2022conjectures} conjectures others. There are more examples from Sun, such as
	$$\sum_{n=1}^\infty 4^{-3n} \frac{(1)_n^3}{n^3(1/2)_n^3} (21 n-8) \left(\frac{43 H_{n-1}^{(3)}}{8}+H_{2 n-1}^{(3)}\right) \stackrel{?}{=} \frac{711 \zeta (5)}{28}-\frac{29 \pi ^2 \zeta (3)}{14}$$
	which \textit{does not follow} from Example \ref{ex5}, and is currently still conjectural. 
	
	\begin{romexample}\label{ex6}
		When $a,b,c,d,e$ are near $0$, we have
		\begin{multline*}\sum_{k\geq 0}\frac{(a+e+2 k+2) (a+1)_k (b+1)_k (c+1)_k (d+1)_k}{(e+1)_{k+1} (a-b+e+1)_{k+1} (a-c+e+1)_{k+1} (a-d+e+1)_{k+1}}= \\
			-\sum_{n\geq 1}(-1)^n P \frac{\splitfrac{(a+1)_{n-1} (b+1)_{n-1}(c+1)_{n-1} (d+1)_{n-1} (-b+e+1)_{n-1} (-c+e+1)_{n-1} }{(-d+e+1)_{n-1} (a-b-c+e+1)_{n-1}(a-b-d+e+1)_{n-1} (a-c-d+e+1)_{n-1}}}{(e+1)_{2 n}(a-b+e+1)_{2 n}(a-c+e+1)_{2 n} (a-d+e+1)_{2 n} (a-b-c-d+2 e+1)_{2 n}}\end{multline*}
	\end{romexample}
	here $P\in \mathbb{Z}[a,b,c,d,e,n]$ is a very long polynomial.\footnote{readers can consult the accompanying file to see the full formula.}
	\begin{proof}
		This is the formula induced by $\textsf{Dougall5F4}(a - e - n, b - e - n, c - e - n, d - e - n, 1 + k + e + 2 n)$, boundary term $\lim_{n\to\infty}\sum_{k\geq 0}F(n,k) = 0$ by Proposition \ref{vanishingboundary}. The formula is $\sum_{k\geq 0}F(0,k) = \sum_{n\geq 1}G(n-1,0)$.
	\end{proof}
	
	Above formula is equivalent to Theorem 31 in Chu \cite{chu2014accelerating}. Compare coefficient on Example \ref{ex6} at the point $(a,b,c,d)=(0,0,0,0)$. LHS will be homogeneous CMZVs of level 1, RHS will be infinite series with hypergeometric part $2^{-10n} \frac{(1)_n^5}{(1/2)_n^5}$ and kind of harmonic numbers involved are $H_n^{(r)}$ and $H_{2n}^{(r)}$. Taking certain linear combinations of above formula, we obtain
	
	\begin{corollary}
		\begin{align*}
			\sum_{n\geq 1}2^{-10n} \frac{(1)_n^5}{n^5(1/2)_n^5}(205 n^2-160 n+32) &= -2\zeta(3) \\
			\sum_{n\geq 1}2^{-10n} \frac{(1)_n^5}{n^5(1/2)_n^5}\left((205 n^2-160 n+32) (H_{2 n-1}-H_{n-1})-41 n+16\right) &= \frac{\pi^4}{60}\\
			\sum_{n\geq 1}2^{-10n} \frac{(1)_n^5}{n^5(1/2)_n^5}\left((205 n^2-160 n+32)(4 H_{2 n-1}^{(2)}-12 H_{n-1}^{(2)})-43 \right) &= -8\zeta(5)\\
			\sum_{n\geq 1}2^{-10n} \frac{(1)_n^5}{n^5(1/2)_n^5}\left(8 (205 n^2-160 n+32)(H_{n-1}^{(3)}+H_{2 n-1}^{(3)})+\frac{125}{n}\right) &= 16\zeta(3)^2\\
			\sum_{n\geq 1}2^{-10n} \frac{(1)_n^5}{n^5(1/2)_n^5}\left(16 (205 n^2-160 n+32)(3H_{n-1}^{(4)}+H_{2 n-1}^{(4)})-\frac{195}{n^2}\right) &=  32\zeta(7)
		\end{align*}
	\end{corollary}
	
	The first equality appeared in \cite{amdeberhan1998hypergeometric}, inspired by it, Sun conjectured remaining four. All of them are essentially proved by Wei in \cite{wei2023some}, although only first three formulas appear explicitly therein, his method works for others as well.  \\[0.04in]
	
	Next two examples (\ref{ex7}, \ref{ex8}) and their corollaries are the most computational intensive part this article. 
	
	\begin{romexample}\label{ex7}
		When $a,b,c,d,e$ are near $0$, \begin{multline*}2\sum_{k\geq 0} \dfrac{(a+3 e+2 k+2) (2 b+2 e+2 k+1) (a+2 e+1)_k (b+e+\frac{1}{2})_k (c+e+1)_k (d+e+1)_k}{(e+1)_{k+1} (a-b+2 e+\frac{1}{2})_{k+1} (a-c+2 e+1)_{k+1} (a-d+2 e+1)_{k+1}} \\ = \sum_{n\geq 1} \frac{\splitfrac{(a+2 e+1)_{2 n-2} (b+e+\frac{1}{2})_n (c+e+1)_{n-1} (d+e+1)_{n-1} (a-b-c+e+\frac{1}{2})_{n-1}}{ (a-b-d+e+\frac{1}{2})_{n-1} (a-c-d+e+1)_{n-1} \times P}}{(e+1)_n (a-b+2 e+\frac{1}{2})_{2 n}(a-c+2 e+1)_{2 n} (a-d+2 e+1)_{2 n} (a-b-c-d+e+\frac{1}{2})_n}\end{multline*}
		here $P\in \mathbb{Z}[a,b,c,d,e,n]$ is a long polynomial\footnote{readers can consult the accompanying file to see the full formula.}.
	\end{romexample}
	\begin{proof}
		This is the formula induced from $F(n,k) = \textsf{Dougall5F4}(a + e + n, 1/2 + b, c, d, 1 + e + k + n)$, the boundary term $\lim_{n\to \infty}\sum_{k\geq 0}F(n,k) = 0$. The formula is $\sum_{k\geq 0}F(0,k) = \sum_{n\geq 1}G(n-1,0)$.
	\end{proof}
	
	When comparing coefficient at $(a,b,c,d,e) = (0,0,0,0,0)$, LHS is (inhomogenous) level 2 CMZV, RHS has form
	\begin{equation}\label{alternatinglevel2form}\sum_{n\geq 1} \frac{1}{\binom{4n}{2n}} \mathbb{Q}[\frac{1}{n},\frac{1}{2n-1}, \hbar_n^{(r)}(1), \hbar_n^{(r)}(1/2)]\end{equation}
	
	Now consider the conjecture by Sun, 
	$$ \sum_{n\geq 1}\frac{1}{n^2 \binom{2 n}{n}} (12 \sum _{j=1}^{n} \frac{(-1)^j}{j^2}-\frac{(-1)^n}{n^2}) = -\frac{11 \pi ^4}{180}$$
	using $$\sum_{j=1}^{2n} \frac{(-1)^j}{j^s} = 2^{1-s} H_{n}^{(s)}-H_{2n}^{(s)} \qquad \sum_{j=1}^{2n-1} \frac{(-1)^j}{j^s} = 2^{1-s} H_{n}^{(s)}-H_{2n}^{(s)} - \frac{1}{(2n)^s} $$ we can split the sum into odd and even parts and then recombine:
	$$\sum _{n\geq 1} \frac{1}{(2n-1)^2 \binom{2(2n-1)}{2n-1}}(6 H_{n}^{(2)}-12H_{2n}^{(2)}-\frac{1}{(2n)^2}+\frac{1}{(2n-1)^2}) + \sum_{n\geq 1} \frac{1}{(2n)^2 \binom{4 n}{2n}}(6 H_n^{(2)}-12H_{2n}^{(2)}-\frac{1}{(2n)^2})$$
	note that $\binom{2(2n-1)}{2n-1} = \frac{n}{4n-1} \binom{4n}{2n}$, we see above sum matches the form in \ref{alternatinglevel2form}. Therefore we hope it is derivable from comparing coefficient of Example \ref{ex7}. This, along with many others are indeed so:
	
	\begin{corollary}[Conjectures in \cite{sun2010conjectures}]
		\begin{align*}\sum _{n=1}^{\infty } \frac{1}{n^2 \binom{2 n}{n}} (12 \sum _{j=1}^{n} \frac{(-1)^j}{j^2}-\frac{(-1)^n}{n^2}) &= -\frac{11 \pi ^4}{180}\\
			\sum _{n=1}^{\infty } \frac{1}{n^2 \binom{2 n}{n}} (24 \sum _{j=1}^{n} \frac{(-1)^j}{j^3}-\frac{17 (-1)^n}{n^3}) & =7 \zeta (5)-\pi ^2 \zeta (3) \\
			\sum _{n=1}^{\infty } \frac{(-1)^n}{n^3 \binom{2 n}{n}}  (10 \sum _{j=1}^n \frac{(-1)^j}{j^2}-\frac{(-1)^n}{n^2})&=\frac{29 \zeta (5)}{6}-\frac{\pi ^2 \zeta (3)}{18} \\
			\sum _{n=1}^{\infty } \frac{1}{n^4 \binom{2 n}{n}} (72 \sum _{j=1}^n \frac{(-1)^j}{j^2}-\frac{(-1)^n}{n^2})&=-\frac{34 \zeta (3)^2}{5}-\frac{31 \pi ^6}{1134} \\
			\sum _{n=1}^{\infty } \frac{1}{n^2 \binom{2 n}{n}} (8 \sum _{j=1}^n \frac{(-1)^j}{j^4}+\frac{(-1)^n}{n^4})&=-\frac{22 \zeta (3)^2}{15}-\frac{97 \pi ^6}{34020}\\
			\sum _{n=1}^{\infty } \frac{(-1)^n}{n^3 \binom{2 n}{n}} (40 \sum _{j=1}^{n} \frac{(-1)^j}{j^3}-\frac{47 (-1)^n}{n^3})&=-\frac{367 \pi ^6}{27216}+6 \zeta (3)^2 \\
			\sum _{n=1}^{\infty } \frac{(-1)^n}{n^3 \binom{2 n}{n}} (110 \sum _{j=1}^n \frac{(-1)^j}{j^4}+\frac{29 (-1)^n}{n^4})&=\frac{221 \pi ^4 \zeta (3)}{180}+\frac{223 \zeta (7)}{24}-\frac{301 \pi ^2 \zeta (5)}{36}\end{align*}
	\end{corollary}
	
	\begin{romexample}
		\label{ex8}
		When $a,b,c,d,e$ are near $0$, $$2^{9}\sum_{k\geq 0}\dfrac{(a+2 e+2 k+1) (a+1)_k (d+\frac{1}{2})_k (b+e+\frac{3}{4})_k (c+e+\frac{1}{4})_k}{(2 e+1)_{k} (a-b+e+\frac{1}{4})_{k+1} \left(a-c+e+\frac{3}{4}\right)_{k+1} \left(a-d+2 e+\frac{1}{2}\right)_{k+1}} = \sum_{n\geq 0} P\times \dfrac{f}{g}$$
		here \begin{multline*}f = (-b+e+\frac{1}{4})_n (b+e+\frac{3}{4})_n (-c+e+\frac{3}{4})_n (c+e+\frac{1}{4})_n(a-b-d+e+\frac{3}{4})_n (a-c-d+e+\frac{5}{4})_{n+1} (-d+2 e+\frac{1}{2})_{2 n} \\ g = (a-b+e+\frac{1}{4})_{n+1} (a-c+e+\frac{3}{4})_{n+1} (2 e+1)_{2 n+1} \left(a-d+2 e+\frac{1}{2}\right)_{2 n+2} (a-b-c-d+2 e+\frac{1}{2})_{2 n+2}\end{multline*}
		and $P\in \mathbb{Z}[a,b,c,d,e,n]$ is a long polynomial\footnote{readers can consult the accompanying file to see the full formula.}.
	\end{romexample}
	\begin{proof}
		This is the formula induced from $F(n,k) = \textsf{Dougall5F4}(a-2 e-2 n+1,b-e-n+\frac{3}{4},c-e-n+\frac{1}{4},d-2 e-2 n+\frac{1}{2},2 e+k+2 n)$, the boundary term $\lim_{n\to \infty}\sum_{k\geq 0}F(n,k) = 0$. The formula is $\sum_{k\geq 0}F(0,k) = \sum_{n\geq 0}G(n,0)$.
	\end{proof}
	
	When expanded around $(a,b,c,d,e) = (0,0,0,0,0)$, LHS is (inhomogenous) CMZV of level 4, RHS is of form 
	$$\sum_{n\geq 0} \frac{\binom{4n}{2n}}{16^{2n}} \mathbb{Q}[\frac{1}{n},\frac{1}{2n+1},\frac{1}{4n+1},\frac{1}{4n+3},\hbar_n^{(r)}(1),\hbar_n^{(r)}(1/2),\hbar_n^{(r)}(1/4),\hbar_n^{(r)}(3/4)]$$

	\begin{corollary}
		\begin{align*}\sum _{n=0}^{\infty } \frac{\binom{2 n}{n}}{(2 n+1) 16^n} (12 \sum _{j=0}^n \frac{(-1)^j}{(2 j+1)^2}-\frac{(-1)^n}{(2 n+1)^2})&=4 \pi  G \\
			\sum _{n=0}^{\infty } \frac{\binom{2 n}{n} }{(2 n+1) 16^n} (24 \sum _{j=0}^{n} \frac{(-1)^j}{(2 j+1)^3}-\frac{17 (-1)^n}{(2 n+1)^3})&=\frac{\pi ^4}{12} \\
			\sum _{n=0}^{\infty } \frac{\binom{2 n}{n}}{(2 n+1)^3 16^n}(9 H_{2 n+1}+\frac{32}{2 n+1})&=40 \beta(4)+\frac{5 \pi  \zeta (3)}{12}\\
			\sum _{n=0}^{\infty } \frac{\binom{2 n}{n} }{(2 n+1)^2 (-16)^n}(10 \sum _{j=0}^n \frac{(-1)^j}{(2 j+1)^2}-\frac{(-1)^n}{(2 n+1)^2})&=G \pi ^2-\frac{\pi  \zeta (3)}{24}\\
			\sum _{n=0}^{\infty } \frac{\binom{2 n}{n} }{(2 n+1)^2 (-16)^n} (40 \sum _{j=0}^{n} \frac{(-1)^j}{(2 j+1)^3}-\frac{47 (-1)^n}{(2 n+1)^3})&=-\frac{85 \pi ^5}{3456}\\
			\sum _{n=0}^{\infty } \frac{\binom{2 n}{n} }{(2 n+1) 16^n} (8 \sum _{j=0}^n \frac{(-1)^j}{(2 j+1)^4}+\frac{(-1)^n}{(2 n+1)^4})&=\frac{11 \pi ^2 \zeta (3)}{120}+\frac{8 \pi  \beta(4)}{3}\\
			\sum _{n=0}^{\infty } \frac{\binom{2 n}{n}}{(2 n+1)^3 16^n} (72 \sum _{j=0}^n \frac{(-1)^j}{(2 j+1)^2}-\frac{(-1)^n}{(2 n+1)^2})&=\frac{7 \pi ^3 G}{3}+\frac{17 \pi ^2 \zeta (3)}{40}\\
			\sum _{n=0}^{\infty } \frac{\binom{2 n}{n} }{(2 n+1)^2 (-16)^n}(110 \sum _{j=0}^n \frac{(-1)^j}{(2 j+1)^4}+\frac{29 (-1)^n}{(2 n+1)^4})&=\frac{91 \pi ^3 \zeta (3)}{96}+11 \pi ^2 \beta (4)-\frac{301 \pi  \zeta (5)}{192}\end{align*}
	\end{corollary}
	
	\begin{proof}
		For each of them, split the summand into even and odd part, then recombine them, note that $$\sum_{j=0}^{2n-1} \frac{(-1)^j}{(2j+1)^s} = \sum_{j=0}^{n} \frac{1}{(4j+1)^s} - \frac{1}{(4j+3)^s} = 4^{-s}(\hbar_{n+1}^{(s)}(1/4) -  \hbar_{n+1}^{(s)}(3/4))$$
		the combined sum can be obtained by comparing coefficients from Example \ref{ex8}.
	\end{proof}
	
	Charlton and Gangl \cite{charlton2023two} proved, via functional equation of single-valued polylogarithm, the third equality above. All others were conjectures due to Sun \cite{sun2010conjectures}, except the fourth, for which we have a stronger form:
	
	\begin{corollary}[Conjecture in \cite{sun2010conjectures}]
		$$\sum _{n=0}^{\infty } \frac{\binom{2 n}{n} }{(2 n+1)^2 (-16)^n}\sum _{j=0}^n \frac{(-1)^j}{(2 j+1)^2}= \frac{\pi^2 G}{10} + \frac{\pi \zeta(3)}{240} + \frac{27\sqrt{3}}{640}L_{-3}(4)$$
	\end{corollary}
	\begin{proof}
		In view of the fourth equality of above corollary, it suffices to prove $$\sum_{n\geq 0}\frac{\binom{2 n}{n}}{(2 n+1)^4 16^n} = \frac{27}{64} \sqrt{3} L_{-3}(4)+\frac{\pi  \zeta (3)}{12}$$
		This has already been proved by Zucker \cite{zucker1985series}. 
	\end{proof}
	
	We note two interesting patterns of above examples. Firstly, when we juxtapose expressions of above two corollaries, one immediate note the coefficients are same:
	$$\sum _{n=0}^{\infty } \frac{\binom{2 n}{n} }{(2 n+1)^2 (-16)^n} (40 \sum _{j=0}^{n} \frac{(-1)^j}{(2 j+1)^3}-\frac{47 (-1)^n}{(2 n+1)^3})=-\frac{85 \pi ^5}{3456}$$
	$$\sum _{n=1}^{\infty } \frac{(-1)^n}{n^3 \binom{2 n}{n}} (40 \sum _{j=1}^{n} \frac{(-1)^j}{j^3}-\frac{47 (-1)^n}{n^3})=-\frac{367 \pi ^6}{27216}+6 \zeta (3)^2$$
	here we have both $40, -47$. The same thing happens for other pairs. Secondly, for every (except the third) entry in previous corollary, the result is a multiple of $\pi$. \par
	
	An explanation of above two patterns might assist us in discovering more similar formulas. \\[0.05in]
	
	Next we give two examples such that the boundary term $\lim_{n\to\infty} \sum_{k\geq 0} F(n,k)\neq 0$ in Proposition $\ref{WZsumformula}$. 
	
	\begin{romexample}\label{ex9}
		For $a,b,c,d$ near $0$, we have
		\begin{multline*}\sum_{k\geq 0} \frac{ (c-d+1)_k (2 a+2 c-d+1)_k (2 b+2 c-d+1)_k}{(d+1)_k (2 c-3 d+1)_{k+1}(a+b+2 c-d+\frac{1}{2})_{k+1}} \\ =  d\sum_{n\geq 1} \frac{\splitfrac{P(-1)^n (1-d)_{n-1} (-a-d+1)_{n-1} (-b-d+1)_{n-1}}{ (c-d+1)_{n-1}(2 a+2 c-d+1)_{n-1} (2 b+2 c-d+1)_{n-1}}}{4 (2 c-3 d+1)_{3 n}(a+c-d+\frac{1}{2})_n (b+c-d+\frac{1}{2})_n (a+b+2 c-d+\frac{1}{2})_n} \\ + \dfrac{\splitdfrac{2^{2 a+2 b+2 c-1} \Gamma (d+1) \Gamma (2 c-3 d+1) \Gamma (-a-b-c+\frac{1}{2})}{ \Gamma (a+c-d+\frac{1}{2}) \Gamma (b+c-d+\frac{1}{2}) \Gamma (a+b+2 c-d+\frac{1}{2})}}{\Gamma (-a-d+1) \Gamma (-b-d+1) \Gamma (c-d+1) \Gamma (2 a+2 c-d+1) \Gamma (2 b+2 c-d+1)}
		\end{multline*}
		here \small{$P = -4 a b d+4 a b n+8 a c^2-28 a c d+28 a c n-4 a c+20 a d^2-40 a d n+6 a d+20 a n^2-6 a n+8 b c^2-28 b c d+28 b c n-4 b c+20 b d^2-40 b d n+6 b d+20 b n^2-6 b n+16 c^3-64 c^2 d+64 c^2 n-12 c^2+76 c d^2-152 c d n+30 c d+76 c n^2-30 c n+2 c-28 d^3+84 d^2 n-18 d^2-84 d n^2+36 d n-3 d+28 n^3-18 n^2+3 n$}
	\end{romexample}
	\begin{proof}
		Let $F(n,k) = \textsf{Watson3F2}(a + d - n, b + d - n, -c + 2 d - 2 n, 1 + 2 c - 3 d + k + 3 n)$, then in accord with Problem \ref{summation_induces_WZ}, one find its WZ-mate $G(n,k)$. It's easy to see all three prerequisites of Proposition \ref{WZsumformula} are satisfied, so we have $$\sum_{k\geq 0}F(0,k) = \sum_{n\geq 1}G(n-1,0) + \lim_{n\to\infty} \sum_{k\geq 0} F(n,k)$$
		let $A = \sum_{k\geq 0}\cdots$, $B = d\sum_{n\geq 1}\cdots$ be the term appearing in statement of the example, then after some computations, one arrives at
		\begin{multline*}A = B +  \dfrac{\Gamma (d+1) \Gamma (a+d+1) \Gamma (b+d+1) \Gamma (2 c-3 d+1) \Gamma (a+b+2 c-d+\frac{1}{2})}{\splitdfrac{(a+d) (b+d) \Gamma (c-d+1) \Gamma (2 a+2 c-d+1) \Gamma (-a-c+d+\frac{1}{2}) }{\Gamma (2 b+2 c-d+1) \Gamma (-b-c+d+\frac{1}{2})}} \\ \times \lim_{n\to \infty}\sum_{k\geq 0} \frac{\splitfrac{\Gamma (-a-c+d-n+\frac{1}{2}) \Gamma (-b-c+d-n+\frac{1}{2})}{\Gamma (c-d+k+n+1) \Gamma (2 a+2 c-d+k+n+1) \Gamma (2 b+2 c-d+k+n+1)}}{\Gamma (a+d-n) \Gamma (b+d-n) \Gamma (d+k-n+1) \Gamma (2 c-3 d+k+3 n+2) \Gamma \left(a+b+2 c-d+k+n+\frac{3}{2}\right)}\end{multline*}
		thus it remains to calculate the last term, this is done in appendix (Lemma \ref{lim_ex9}). 
	\end{proof}
	
	Note that Example \ref{ex9} reduces to Watson's $_3F_2$ formula when $d=0$. Comparing coefficient of $a^0b^0c^0d^1$ on both sides of Example \ref{ex9} produces $$\sum _{k=1}^{\infty } \frac{2 (1)_k H_k }{(-2 k-1) (k+1) \left(\frac{1}{2}\right)_k} = \sum_{n\geq 1} \frac{(-1)^n \left(28 n^3-18 n^2+3 n\right) (1)_n^6}{2 n^6 (\frac{1}{2})_n^3 (3n)!}$$
	LHS easily evaluates to $-7\zeta(3)$ (more generally, such sum are CMZV of level 2 or 4, see\footnote{actually invoking this proposition is not necessary if one develops the this WZ-technique further, we postpone this to a later occasion.} Proposition \ref{binomCMZVs}), so we prove the first of the following corollary.
	\begin{corollary}[Conjectures in \cite{sun2013products}, \cite{sun2021book}, \cite{sun2010conjectures}]
		Let $a_n = \dfrac{(-\frac{1}{27})^n (1)_n^5}{n^5 (\frac{1}{3})_n (\frac{1}{2})_n^3 (\frac{2}{3})_n}$, we have \begin{align*}\sum_{n\geq 1} a_n (28 n^2-18 n+3) &= -14\zeta(3) \\
			\sum_{n\geq 1} a_n \left((28 n^2-18 n+3) \left(4 H_{2 n-1}-3 H_{n-1}\right)-20 n+6 \right) &= \frac{\pi^4}{2} \\
			\sum_{n\geq 1} a_n \left((28 n^2-18 n+3) (2 H_{2 n-1}^{(2)}-3 H_{n-1}^{(2)})-2\right) &= -31 \zeta (5) \\
			\sum_{n\geq 1} a_n \left((28 n^2-18 n+3) (3 H_{n-1}^{(3)}+4 H_{2 n-1}^{(3)})+\frac{4}{n}\right) &= -49 \zeta (3)^2
		\end{align*}
	\end{corollary}
	\begin{proof}
		All follow from comparing coefficient of Example \ref{ex9}. 
	\end{proof}

	\begin{romexample}\label{ex10}For $a,b$ near $0$, we have
		\begin{multline*}\sum_{k\geq 0} \frac{2 \left(-\frac{1}{3}\right)^k \left(-a+b+\frac{1}{2}\right)_k (2 a+b+1)_k}{(4 a+2 b+2 k+1) (b+1)_k \left(2 a+b+\frac{1}{2}\right)_k} = \frac{4^a 3^{-a+b-1/2} \Gamma (b+1) \Gamma \left(2 a+b+\frac{1}{2}\right) \Gamma \left(a-b+\frac{1}{2}\right)}{\Gamma (a+1) \Gamma (2 a+b+1)} \\ - b\sum_{n\geq 1} \frac{3^n \left(a+\frac{1}{2}\right)_n (a+1)_n^2 (2 a+b+1)_{2 n}}{(a+n) (2 a+b+2 n-1) (2 a+1)_{2 n} \left(a-b+\frac{1}{2}\right)_n \left(2 a+b+\frac{1}{2}\right)_{2 n}}
		\end{multline*}
	\end{romexample}
	\begin{proof}
		Like last example, the gamma product originates again from $\lim_{n\to\infty} \sum_{k\geq 0} F(n,k)$. Let $$F(n,k) = \frac{(-1)^k 3^{n-k} \Gamma (a+n+1)^2 \Gamma \left(-a+b+k-n+\frac{1}{2}\right) \Gamma (2 a+b+k+2 n+1)}{\Gamma \left(-a-n+\frac{1}{2}\right) \Gamma (2 a+2 n+1) \Gamma (b+k+1) \Gamma \left(2 a+b+k+2 n+\frac{3}{2}\right)}$$
		one checks $$G(n,k) = \frac{3 (b+k) (2 a+b+k+2 n+2)}{(2 a-2 b-2 k+2 n+1) (4 a+2 b+2 k+4 n+3)} F(n,k)$$ is its WZ-mate. All three prerequisites of Proposition \ref{WZsumformula} are met. 
		We will calculate $$\lim_{n\to\infty} \sum_{k\geq 0} F(n,k) = \sqrt{\pi } 3^{-a+b-\frac{1}{2}} \cos (\pi a) \sec (\pi  (-a+b))$$ in Lemma \ref{lim_ex10} in appendix. Dividing both sides $$\sum_{k\geq 0}F(0,k) = \sum_{n\geq 1} G(n-1,0) + \lim_{n\to\infty} \sum_{k\geq 0} F(n,k)$$ by $\frac{\Gamma (a+1)^2 \Gamma \left(-a+b+\frac{1}{2}\right) \Gamma (2 a+b+1)}{\Gamma \left(\frac{1}{2}-a\right) \Gamma (2 a+1) \Gamma (b+1) \Gamma \left(2 a+b+\frac{1}{2}\right)}$ gives the formula. 
	\end{proof}

	Comparing coefficient of $a^0b^1$ on both sides, we have
	$$\sum _{k=0}^{\infty } -\frac{4 \left(-\frac{1}{3}\right)^k}{(2 k+1)^2} = \frac{\pi  \log 3}{\sqrt{3}} - \sum _{n=1}^{\infty } \frac{3^n (1)_n^2}{n (2 n-1) \left(\frac{1}{2}\right)_{2 n}}$$
	LHS is 
	$$2\sqrt{3}i\sum\limits_{k = 1}^\infty  {\frac{{{{(\frac{i}{{\sqrt 3 }})}^k}}}{{{k^2}}}[1 - {{( - 1)}^k}]}  = 2\sqrt{3}i(\Li_2(\frac{i}{{\sqrt 3 }}) - \Li_2(\frac{{ - i}}{{\sqrt 3 }}))$$
	which can be easily shown\footnote{There are many ways. For example: using an \textit{ad hoc} manipulation of functional equation of $\Li_2$; or invoking the fact that $\Li_n(i/\sqrt{3})$ is CMZV of level 6 weight $n$, so they can be expressed using our chosen CMZV basis} $=\frac{\pi  \log (3)}{\sqrt{3}}-\frac{15 L_{-3}(2)}{2}$, hence \begin{corollary} [Conjecture in \cite{sun2021book}, \cite{sun2010conjectures}]
		$$L_{-3}(2)=\frac2{15}\sum _{k=1}^\infty\frac{48^k}{k(2k-1)\binom{4k}{2k}\binom{2k}k}$$
	\end{corollary}
	
	Actually, our choice of $F(n,k)$ in proof of Example \ref{ex10} is, like all other examples in this article, also inspired from a hypergeometric summation formula. This formula can be recovered by setting $b=0$:
	$$\frac{2}{1+4a}\times  \pFq{2}{1}{1/2-a,1+2a}{3/2+2a}{-\frac{1}{3}} = \frac{3^{-a-\frac{1}{2}} 4^a \Gamma \left(a+1/2\right) \Gamma \left(2 a+\frac{1}{2}\right)}{\Gamma (a+1) \Gamma (2 a+1)}$$
	discovered by Goursat in 1881 \cite{goursat1881equation}. 

	\newpage 
	
	\appendix
	\section{Proof of Theorem \ref{CMZVsum1}}
	We present here the technical proof of Theorem \ref{CMZVsum1}, which we give due to a lack of suitable reference. \par
	
	For convenience, we restate the statement to be proved:
	\begin{theorem}
		Let $0<\gamma_0\leq 1$, $K_i$ positive integers, $\vec{\gamma_i}$ and $\vec{a_i}$ rational vectors such that each component is $0< \cdot \leq 1$. Let $N$ be an integer such that all $N\vec{\gamma_i}, N\vec{a_i}, N\gamma_0, Na_0, N/K_i$ are in $\mathbb{Z}$. The following series, regularized, is a level $N$ (possibly inhomogeneous) CMZV
		$$\sum_{n_0>0} \frac{e^{2\pi i  a_0 n_0}}{(n_0+\gamma_0)^{s_0}} H(\vec{s_1},\vec{\gamma_1},\vec{a_1};K_1 n_0) \cdots H(\vec{s_m},\vec{\gamma_m},\vec{a_m};K_m n_0) \in \sum_{1\leq i\leq w} \CMZV{N}{i}$$
		where $w = |\vec{s_1}| + \cdots + |\vec{s_i}| + s$. If $\gamma_0 = 1/N$, then it is a homogenous CMZV of weight $w$. 
	\end{theorem}
	
	We tacitly assume all entries of $\vec{s}$ are in $\in \mathbb{N}$, all entries of $\vec{\gamma}$ are $\in \mathbb{Q}\cap (0,1]$. 
	\begin{lemma}
		Given positive integer $K$, $H(\vec{s},\vec{\gamma},\vec{a};Kn)$ can be expressed in terms of finitely linear combinations of various $H(\vec{s},\vec{\gamma}',K\vec{a};n)$. 
	\end{lemma}
	\begin{proof}
		This is straightforward. By definition, $$\begin{aligned} H(\vec{s},\vec{\gamma},\vec{a};Kn) &= \sum_{Kn\geq n_1+\gamma_1 > n_2+\gamma_2 > \cdots > n_k+\gamma_k > 0} \frac{e(a_1 n_1)\cdots e(a_k n_k)}{(n_ 1+\gamma_1)^{s_ 1}\cdots (n_k+\gamma_k)^{s_k}} \\
			&= \sum_{\vec{r} \in [0,K-1]^k} \quad \sum_{Kn\geq Km_1 + r_1+\gamma_1 > Km_2+r_2+\gamma_2 > \cdots > Km_k+r_k+\gamma_k > 0} \prod_i \frac{e(a_i(Km_i+r_i))}{(Km_i+r_i+\gamma_i)^{s_i}} \\
			&= \sum_{\vec{r} \in [0,K-1]^k} \quad (\prod_i e(a_i r_i) K^{-s_i}) \sum_{n\geq m_1 + (r_1+\gamma_1)/K > \cdots > m_k+(r_k+\gamma_k)/K > 0} \prod_i \frac{e(Ka_i m_i)}{(m_i+(r_i+\gamma_i)/K)^{s_i}} \\
			&= \sum_{\vec{r} \in [0,K-1]^k} \quad (\prod_i e(a_i r_i) K^{-s_i}) H(\vec{s}, \frac{\vec{\gamma}+\vec{r}}{K}, K\vec{a};n)
		\end{aligned}$$
		by our assumption that each entries of $\vec{\gamma}$ is in $(0,1]$ implies the same is true for entries of $(\vec{\gamma} + \vec{r})/K$. 
	\end{proof}
	
	\begin{lemma}
		The product of two $H(\vec{s_1},\vec{\gamma_1},\vec{a_1};n) H(\vec{s_2},\vec{\gamma_2},\vec{a_2};n)$ can be expressed in terms of finitely $\mathbb{Z}$-linear combinations of various $H(\vec{s}',\vec{\gamma}',\vec{a}';n)$.
	\end{lemma}
	\begin{proof}
		This is essentially stuffle product. Let $\mathfrak{A}$ be the $\mathbb{Q}$-algebra generated by alphabets $$\mathcal{A} = \{w_{s,\gamma,a} | s\in \mathbb{N},\gamma \in \mathbb{Q}\cap (0,1], a\in \mathbb{C} \}.$$ Define an operation $\ast$ on $\mathfrak{A}$, which we stipulate to be $\mathbb{Q}$-linear, inductively as follows:
		$$(u_1 v_1) \ast (u_2 v_2) = u_1 (v_1 \ast u_2 v_2) + u_2 (u_1 v_1 \ast v_2) + f(u_1,u_2) (v_1\ast v_2) $$
		for $u_1,u_2 \in \mathcal{A}$ and
		$$f(w_{s_1,\gamma_1,a_1},w_{s_2,\gamma_2,a_2}) = \begin{cases}0 & \text{ if }\gamma_1\neq \gamma_2 \\ w_{s_1+s_2,\gamma_1,a_1+a_2} & \text{ if }\gamma_1 = \gamma_2\end{cases}$$
		The operation is easily seen to be associative. Furthermore, for each $n \in \mathbb{N}$ and $v\in \mathfrak{A}$, we can define $H(v,n)$ via
		$$H(w_{s_1,\gamma_1,a_1}\cdots w_{s_k,\gamma_k,a_k};n) = H(\vec{s},\vec{\gamma},\vec{a};n)$$
		and distribute this linearly to all of $\mathfrak{A}$. 
		Using familiar arguments from stuffle product, one sees this map is a homomorphism, i.e. $H(v_1;n)H(v_2;n) = H(v_1\ast v_2;n)$. This is the claim. 
	\end{proof}
	
	\begin{lemma}
		We have the following asymptotic expansion:
		$$H(\vec{s},\vec{\gamma},\vec{a};M) = \sum_{i>0} c_i (\log M + \gamma)^i + c_0 + o(1) \qquad M\to \infty$$
		with $c_i \in \CMZV{N}{w-i}$, almost all $c_i = 0$. Here $N$ is an integer such that $N \vec{\gamma},N\vec{a}$ are integral, $w = \sum_i s_i$. 
	\end{lemma}
	\begin{proof}
		Let $R$ be the least common denominator of entries of $\gamma$. Using $m_i =R(n_i+\gamma_i)$, we have
		$$\begin{aligned}H(\vec{s},\vec{\gamma},\vec{a};M) &= \sum_{\substack{n_1+\gamma_1 > n_2+\gamma_2 > \cdots > n_k+\gamma_k > 0 \\ n_1<M}} \prod_i \frac{e(a_i n_i)}{(n_ i+\gamma_i)^{s_i}} \\
			&= R^{s_1+\cdots+s_k} \sum_{\substack{m_1>m_2\cdots > m_k>0 \\ m_1/R-\gamma_1<M \\ m_i\equiv R\gamma_i \pmod{R}}} \prod_i \frac{e(a_i (\frac{m_i}{R}-\gamma_i))}{m_i^{s_i}}\end{aligned}$$
		the characteristic function of $m_i \equiv R\gamma_i \pmod{R}$ is $\frac{1}{R}\sum_{x_i \pmod{R}} e(x_i \frac{m_i-R\gamma_i}{R})$, so above is
		$$e(-\vec{a}\cdot \vec{\gamma}) R^{s_1+\cdots+s_k-k} \sum_{\vec{x} \in [0,R-1]^k} e(-\vec{x}
		\cdot \vec{\gamma}) \sum_{\substack{m_1>m_2\cdots > m_k>0 \\ m_1/R-\gamma_1<M}} \prod_i \frac{e(m_i (x_i+a_i)/R)}{m_i^{s_i}}$$
		the innermost sum is partial sum of multiple polylogarithm
		$$\sum_{m_1>m_2>\cdots m_k > 0} \prod_i \frac{b_i^{m_i}}{m_i^{s_i}}$$
		with $b_i = e(\frac{x_i+a_i}{R})$ which is an $N$-th root of unity, it's known, from stuffle regularization of MZVs, that $$\sum_{\substack{m_1>m_2\cdots > m_k>0 \\ m_1/R-\gamma_1<M}} \prod_i \frac{e(m_i (x_i+a_i)/R)}{m_i^{s_i}}$$ has the same form of asymptotic expansion as in statement of this lemma, with stuffle regulator $\log R \in \CMZV{N}{1}$, this completes the proof.
	\end{proof}
	
	We need to introduce one more definition: for an ordered list $\vec{\sigma}$ of length $k-1$, whose element are either the symbol $\geq$ or $>$, We denote
	$$H(\vec{s},\vec{\gamma},\vec{a},\vec{\sigma};N) = \sum_{N> n_1 \square_1 n_2 \square_2 \cdots \square_{k-1} n_k > 0} \frac{e(a_1 n_1)\cdots e(a_k n_k)}{(n_ 1+\gamma_1)^{s_ 1}\cdots (n_k+\gamma_k)^{s_k}} \qquad \vec{\sigma} = (\square_1,\cdots,\square_{k-1})$$
	Given $\vec{\gamma}$, we can associate a $\vec{\sigma}_\gamma$ as follows: its $i$-th entry is defined to be
	$$(\vec{\sigma}_\gamma)_i = \begin{cases}> &\text{ if }\gamma_k \leq \gamma_{k+1} \\ \geq &\text{ if }\gamma_k > \gamma_{k+1}\end{cases}$$
	note that $H(\vec{s},\vec{\gamma},\vec{a},\vec{\sigma}_\gamma;N) = H(\vec{s},\vec{\gamma},\vec{a};N)$. 
	
	\begin{lemma}
		We have the following asymptotic expansion:
		$$H(\vec{s},\vec{\gamma},\vec{a},\vec{\sigma};M) = \sum_{i>0} c_i (\log M + \gamma)^i + c_0 + o(1) \qquad M\to \infty$$
		with $c_0 \in \sum_{1\leq v\leq w}\CMZV{N}{v}$, almost all $c_i = 0$. Here $N$ is an integer such that $N \vec{\gamma},N\vec{a}$ are integral, $w = \sum_i s_i$. 
	\end{lemma}
	\begin{proof}
		We prove this statement by induction on the quantity
		$$(\text{number of different entries between }\vec{\sigma} \text{ and }\vec{\sigma}_\gamma) + (\text{the depth }k)$$
		We first look at the base case: when the depth $k$ is $1$, $$H(\vec{s},\vec{\gamma},\vec{a},\vec{\sigma}_\gamma;M) = \sum_{M>n_1>0} \frac{e(a_1 n_1)}{(n_1+\gamma_1)^{s_1}}$$
		for which the statement obviously holds; when the number of different entries between $\vec{\sigma}$ and $\vec{\sigma}_\gamma$ is $0$, we have $H(\vec{s},\vec{\gamma},\vec{a},\vec{\sigma}_\gamma;M) = H(\vec{s},\vec{\gamma},\vec{a};M)$, so the statement holds by the previous lemma. \par
		Therefore we assume $\vec{\sigma}\neq \vec{\sigma}_\gamma$, let $j$ be smallest index for which these two are different, \begin{itemize}
			\item If $\gamma_j \leq \gamma_{j+1}$, then $\sigma_j$ is the symbol $\geq$, let $\vec{\sigma}'$ denote the vector by changing this to $>$, also let $\varepsilon = 1$.
			\item If $\gamma_j > \gamma_{j+1}$, then $\sigma_j$ is the symbol $>$, let $\vec{\sigma}'$ denote the vector by changing this to $\geq$, also let $\varepsilon = -1$.
		\end{itemize}
		then (with $\vec{\sigma} = (\square_1,\cdots,\square_{k-1})$):
		$$H(\vec{s},\vec{\gamma},\vec{a},\vec{\sigma};M) = H(\vec{s},\vec{\gamma},\vec{a},\vec{\sigma}';M) + \varepsilon \times \sum_{N> n_1 \square_1 n_2 \cdots n_j = n_{j+1} \square_{j+1} \cdots \square_{k-1} n_k> 0} \frac{e(a_1 n_1)\cdots e(a_k n_k)}{(n_ 1+\gamma_1)^{s_ 1}\cdots (n_k+\gamma_k)^{s_k}} $$
		using partial fraction, we can split $$\frac{1}{(n+\gamma_j)^{s_j} (n+\gamma_{j+1})^{s_{j+1}}}$$
		from which we see that the latter sum in above equation has lower depth, so by induction hypothesis, it has the desired asymptotic expansion; for the first term $H(\vec{s},\vec{\gamma},\vec{a},\vec{\sigma}';M)$, induction assumption also applies, since $j$-th entry of $\vec{\sigma}'$ and $\vec{\sigma}_\gamma$ now agree. This says $H(\vec{s},\vec{\gamma},\vec{a},\vec{\sigma};M)$ has the desired expansion, concluding the proof.
	\end{proof}
	
	Now we can quickly conclude the proof of Theorem \ref{CMZVsum1}. By first two lemmas, it suffices to prove the statement for $$\sum_{N> n_0>0} \frac{e^{2\pi i  a_0 n_0}}{(n_0+\gamma_0)^{s_0}} H(\vec{s},\vec{\gamma},\vec{a},n_0) = H(\vec{s}',\vec{\gamma}',\vec{a}',\vec{\sigma};N)$$
	here $\vec{s}'$ (resp. $\vec{\gamma}',\vec{a}'$) is obtained by prepending $\vec{s}$ (resp. $\vec{\gamma},\vec{a}$) with $s_0$ (resp. $\gamma_0,a_0$); $\vec{\sigma}$ is obtained by prepending $\vec{\sigma}_\gamma$ with symbol $>$. The claim about its regularized value follows immediately from previous lemma. \par Finally, if $\gamma_0 = 1/N$, then it is smaller or equal to all $\gamma_i$ appearing, so $H(\vec{s}',\vec{\gamma}',\vec{a}',\vec{\sigma};N) = H(\vec{s}',\vec{\gamma}',\vec{a}';N)$, it being a weight-homogeneity of CMZV follows from the penultimate lemma. 
	
	\section{Some asymptotics of hypergeometric series}
	Here we calculate the boundary terms $\lim_{n\to\infty}\sum_{k\geq 0} F(n,k)$ in Example \ref{ex9} and \ref{ex10}. We put the proofs here since they have quite different flavours than our main result.
	\begin{lemma}\label{lim_ex9} The limit at the proof of Example \ref{ex9}, 
		$$\lim_{n\to \infty}\sum_{k\geq 0} \frac{\splitfrac{\Gamma (-a-c+d-n+\frac{1}{2}) \Gamma (-b-c+d-n+\frac{1}{2})}{\Gamma (c-d+k+n+1) \Gamma (2 a+2 c-d+k+n+1) \Gamma (2 b+2 c-d+k+n+1)}}{\Gamma (a+d-n) \Gamma (b+d-n) \Gamma (d+k-n+1) \Gamma (2 c-3 d+k+3 n+2) \Gamma \left(a+b+2 c-d+k+n+\frac{3}{2}\right)} $$ when $a,b,c,d$ are near $0$, equals $$2^{2 a+2 b+2 c-1}\Gamma (-a-b-c+\frac{1}{2}) \sin (\pi  (a+d)) \sin (\pi  (b+d))  \sec (\pi  (a+c-d)) \sec (\pi  (b+c-d))$$
	\end{lemma}
	\begin{proof}Let $f(n,k)$ be the term inside the summand. First note that, using same method as in proof of Proposition \ref{vanishingboundary}, $\lim_{n\to\infty} \sum_{k=0}^{n-1} f(n,k) = 0$ (each term $f(n,k)$ decays geometrically fast for $0\leq k \leq n-1$). Therefore we only need to calculate \begin{align*} &\sum_{k=n}^{\infty} f(n,k) = \sum_{k=0}^{\infty} f(n,k+n) \\ & = \frac{\Gamma (c-d+2 n+1) \Gamma \left(-a-c+d-n+\frac{1}{2}\right) \Gamma (2 a+2 c-d+2 n+1) \Gamma \left(-b-c+d-n+\frac{1}{2}\right) \Gamma (2 b+2 c-d+2 n+1)}{\Gamma (d+1) \Gamma (a+d-n) \Gamma (b+d-n) \Gamma (2 c-3 d+4 n+2) \Gamma \left(a+b+2 c-d+2 n+\frac{3}{2}\right)} \\ &\quad \times \underbrace{\sum_{k\geq 0} \frac{(c-d+2 n+1)_k (2 a+2 c-d+2 n+1)_k (2 b+2 c-d+2 n+1)_k}{(d+1)_k (2 c-3 d+4 n+2)_k \left(a+b+2 c-d+2 n+\frac{3}{2}\right)_k}}_S\end{align*}
		Using $\Gamma(x)\Gamma(1-x) = \pi\csc \pi x$, one sees the gamma factor exterior of $S$, as $n\to\infty$, is $$\sim \frac{\sqrt{\pi } \sin (\pi  (a+d)) \sin (\pi  (b+d)) \sec (\pi  (a+c-d)) \sec (\pi  (b+c-d))}{\Gamma (d+1)}2^{a+b-c+4 d-4 n-2} n^{-a-b-c+d} $$
		It remains to find the asymptotic leading term of $S$. We first note the following easy observation: if $A,B,C,D$ are independent of $n,k$ and $A+B=C+D$ then $$\frac{(A+n)_k (B+n)_k}{(C+n)_k (D+n)_k} = 1 + O(\frac{1}{n})$$ with $O$-term \textit{uniform} in $k\geq 0$. Using this, we see
		\begin{align*}S &\sim \sum_{k\geq 0} \frac{(2 b+2 c-d+2 n+1)_k \left(a-b+c-d+2 n+\frac{1}{2}\right)_k}{(d+1)_k (2 c-3 d+4 n+2)_k} \\
			&= \pFq{3}{2}{1,1+2b+2c-d+2n,1/2+a-b+c-d+2n}{1+d,2+2c-3d+4n}{1} \\
			&= \frac{ d\times \Gamma (2 c-3 d+4 n+2)}{\Gamma (-2 b-2 d+2 n+1) \Gamma (2 b+2 c-d+2 n+1)} \\ &\quad \times \int_{[0,1]^2}(1-s)^{d-1} (1-t)^{-2 b-2 d+2 n}t^{2 b+2 c-d+2 n} (1-s t)^{-a+b-c+d-2 n-1/2} dsdt\end{align*}
		where we used integral representation of hypergeometric function. Make substitution $t=z/(1+z), s\mapsto 1-s$ to the integral, it equals $$\int_0^\infty z^{2 b+2 c-d} (z+1)^{a-b-c+2 d-3/2} (\frac{z}{1+z})^{2n}\left( \int_0^1 s^{d-1} (s z+1)^{-a+b-c+d-2 n-1/2} ds \right) dz$$
		from which we see only large $z$ contributes the the asymptotic, so we assume $z$ is large, then 
		\begin{align*}\int_0^1 s^{d-1} (s z+1)^{-a+b-c+d-2 n-1/2} ds &= (nz)^{-d} \int_0^{nz} s^{d-1} (1+\frac{s}{n})^{-1/2-a+b-c+d-2n} ds \\ &\sim (nz)^{-d} \int_0^{
				\infty} s^{d-1} e^{-2s} ds = (2nz)^{-d}\Gamma(d) \end{align*}
		so the double integral is $$\sim (2n)^{-d}\Gamma(d) \int_0^\infty z^{2 b+2 c-d} (z+1)^{a-b-c+2 d-3/2} (\frac{z}{1+z})^{2n} z^{-d} dz$$
		this is Euler's beta function, so we can easily calculate the leading term of $S$, which is $$S \sim \pi^{-1/2}\Gamma (d+1) \Gamma \left(-a-b-c+\frac{1}{2}\right) 2^{a+b+3 c-4 d+4 n+1} n^{a+b+c-d}$$ This concludes the proof. 
	\end{proof}

	\begin{lemma}\label{lim_ex10}
		When $a,b$ are near $0$, 
		$$\lim_{n\to \infty}\sum_{k\geq 0} \frac{(-1)^k 3^{n-k} \Gamma (a+n+1)^2 \Gamma \left(-a+b+k-n+\frac{1}{2}\right) \Gamma (2 a+b+k+2 n+1)}{\Gamma \left(-a-n+\frac{1}{2}\right) \Gamma (2 a+2 n+1) \Gamma (b+k+1) \Gamma \left(2 a+b+k+2 n+\frac{3}{2}\right)}$$ is $$\sqrt{\pi } 3^{-a+b-\frac{1}{2}} \cos (\pi a) \sec (\pi  (-a+b))$$
	\end{lemma}
	\begin{proof}
		We use the same technique previously. The sum equals
		$$\frac{3^n \Gamma (a+n+1)^2 \Gamma \left(-a+b-n+\frac{1}{2}\right) \Gamma (2 a+b+2 n+1)}{\Gamma (b+1) \Gamma \left(-a-n+\frac{1}{2}\right) \Gamma (2 a+2 n+1) \Gamma \left(2 a+b+2 n+\frac{3}{2}\right)} \times \pFq{3}{2}{1,1/2-a+b-n,1+2a+b+2n}{1+b,3/2+2a+b+2n}{-\frac{1}{3}}$$
		which is $$\sim \frac{3^n 2^{-2 a-2 n} \cos(\pi a) n^{b+\frac{1}{2}} \sec(\pi(b-a))}{\Gamma (b)}\times \int_{[0,1]^2} (1-s)^{b-1} \left(1+\frac{s t}{3}\right)^{a-b+n-1/2} t^{2 a+b+2 n} (1-t)^{-1/2}  ds dt$$
		Let $I$ be the double integral, we make substitutions $t=z/(1+z),s\mapsto 1-s$,
		$$I = \int_0^\infty (\frac{z}{1+z})^{2a+b+2n} (\frac{1}{1+z})^{-3/2} \left(\int_0^1 s^{b-1} \left(1+\frac{1-s}{3}\frac{z}{1+z}\right)^{-1/2+a-b+n} ds \right) dt$$
		
		Only large $z$ contributes to the asymptotic, the inner integral is
		\begin{align*}&= n^{-b} \int_0^{n} s^{b-1}\left(1+\frac{1-s/n}{3}\frac{z}{1+z}\right)^{-1/2+a-b+n} ds\\
			&\sim n^{-b}\left(\frac{3+4z}{3(1+z)}\right)^{-1/2+a-b+n}\int_0^\infty s^{b-1} \exp(-\frac{sz}{3+4z}) ds \\
			&= n^{-b} z^{-b} \Gamma (b) 3^{-a+b-n+\frac{1}{2}} (4 z+3)^{a+n-\frac{1}{2}} (z+1)^{-a+b-n+\frac{1}{2}}\end{align*}
		so \begin{align*}I&\sim 3^{-a+b-n+1/2}n^{-b}\Gamma(b) \int_0^\infty \left( \frac{z^2(3+4z)}{(1+z)^3}\right)^n z^a (1+z)^{-2a-1} (4z+3)^{a-1/2} dz \\
			&= 3^{-a+b-n+1/2}n^{1-b}\Gamma(b) \int_0^\infty \left( \frac{n^2 z^2(3+4nz)}{(1+nz)^3}\right)^n (nz)^a (1+nz)^{-2a-1} (4nz+3)^{a-1/2} dz \\
			&\sim 3^{-a+b-n+1/2}n^{1-b}\Gamma(b) \int_0^\infty 4^n \exp(-\frac{9}{4z}) (nz)^a (nz)^{-2a-1} (4nz)^{a-1/2} \\
			&= \sqrt{\pi } 4^{a+n} n^{-b-\frac{1}{2}} \Gamma (b) 3^{-a+b-n-\frac{1}{2}}\end{align*}
		this completes the proof. 
	\end{proof}

	\bibliographystyle{plain} 
	\bibliography{ref.bib} 
	
\end{document}